\documentclass[aps,prb,twocolumn,floatfix,nofootinbib,reprint,superscriptaddress]{revtex4-1}
\usepackage{amsmath,amssymb,color,comment}
\usepackage[makeroom]{cancel}
\usepackage[caption=false]{subfig}
\usepackage{mathrsfs}
\usepackage{graphicx}
\usepackage{subfig}
\usepackage[countmax]{subfloat}
\usepackage[english]{babel}
\usepackage[bookmarks=true,colorlinks,linkcolor=blue,urlcolor=blue,citecolor=blue]{hyperref}
\usepackage[normalem]{ulem}

\usepackage{bm}

\newcommand{\blue}[1]{{\color{blue}#1}}

\begin{document}

\title{Fast logarithmic Fourier-Laplace transform of nonintegrable functions}

\author{Johannes Lang}

\affiliation{Physik Department, Technische Universit\"at M\"unchen, 85747 Garching, Germany}

\author{Bernhard Frank}

\affiliation{Physik Department, Technische Universit\"at M\"unchen, 85747 Garching, Germany}
\affiliation{Max-Planck-Institut f\"ur Physik komplexer Systeme, 01187 Dresden, Germany}

%\author{Wilhelm Zwerger}

%\affiliation{Physik Department, Technische Universit\"at M\"unchen, 85747 Garching, Germany}

\date{\today}

\begin{abstract}
We present an efficient and very flexible numerical fast Fourier-Laplace transform, that extends the logarithmic Fourier transform (LFT) introduced by Haines and Jones [Geophys. J. Int. 92(1):171 (1988)] for functions varying over many scales to nonintegrable functions. In particular, these include cases of the asymptotic form $f(\nu\to0)\sim\nu^a$ and $f(|\nu|\to\infty)\sim\nu^b$ with arbitrary real $a>b$. Furthermore, we prove that the numerical transform converges exponentially fast in the number of data points, provided that the function is analytic in a cone $|\Im{\nu}|<\theta|\Re{\nu}|$ with a finite opening angle $\theta$ around the real axis and satisfies $|f(\nu)f(1/\nu)|<\nu^c$ as $\nu\to 0$ with a positive constant $c$, which is the case for the class of functions with power-law tails. Based on these properties we derive ideal transformation parameters and discuss how the logarithmic Fourier transform can be applied to convolutions. The ability of the logarithmic Fourier transform to perform these operations on multiscale (non-integrable) functions with power-law tails with exponentially small errors makes it the method of choice for many physical applications, which we demonstrate on typical examples. These include benchmarks against known analytical results inaccessible to other numerical methods, as well as physical models near criticality.
\end{abstract}

\maketitle

\section{Introduction}
In physics, one is often confronted with the need to Fourier transform or convolve functions that are either only numerically available or whose exact transformation is not known. Since the reinvention of the fast Fourier transform (FFT) by Cooley and Tukey~\cite{FFT1965}, which reduces the numerical cost for both of these operations from $\mathcal{O}(N^2)$ to $\mathcal{O}(N \log_2{N})$, where $N$ denotes the number of grid points, the FFT has been established as the standard  method for most situations. However, it necessarily requires an equidistant grid, which is quite inconvenient for many applications in theoretical physics. There, one frequently has to deal with slowly (i.e. algebraically) decaying functions, while the opposite limit of small arguments contains a lot of physical information. An example is provided by
Green's functions in many-body problems with short-range interactions~\cite{zwer14varenna}. 
%correlation functions in the context of ultracold Fermions with short-range interactions. The latter give rise to characteristic power laws at frequencies or momenta which exceed all intrinsic energy or inverse lengths scales. For example the Tan energy theorem~\cite{tan08energy} entails the algebraic tail $n_\sigma(k) \sim \mathcal{C}/k^4$ in the momentum distribution of the Fermions of spin~$\sigma$, where the prefactor is given by the observable Tan contact density $\mathcal{C}$. Regarding the phase diagram and transport properties, however, \red{a lot of}\green{most}\blue{the relevant} information is contained in the low-energy spectrum\cite{zwer14varenna}. A quantitative description of the physics therefore has to account for both limits correctly. A similar problem is encountered in the physics of glass transitions, where experimental data is taken over time scales spanning up to eight orders of magnitude~\cite{Kob1995}.\\
To implement an FFT under such circumstances, it is necessary to use a fine grid for small arguments that extends to very high frequencies, which is of course not very practicable due to the huge number of required data points. Consequently, a number of alternative methods have been introduced in the literature: Sometimes, sufficient knowledge about the asymptotic behavior at large arguments can be gained, subtracted and treated separately, such that the remainder of the function under consideration decays fast enough to be amenable to the application of an FFT~\cite{BDMC2011, BDMC2013a, BDMC2013b}. More often, however, it is necessary to waive the advantages of the FFT in favor of a more flexible sampling, specifically adapted to the problem. This, however requires to apply a discrete Fourier transform (DFT) with $\mathcal{O}(N^2)$ numerical complexity~\cite{Num_recipes}.
%Due to this fast growth of computational effort, there is a strong desire to reduce the number of data points: One way to achieve this is to combine the DFT with a spline interpolation~\cite{Haussmann1994}, which can be done without changing the scaling of the numerical effort using arbitrary spline orders \cite{Frank_thesis}.\\

A combination of the best of both worlds, i.e. an $N \log_2{N}$ scaling on a logarithmic grid, which is able to cover all physically relevant orders of magnitude, has first been proposed by Haines and Jones in form of the logarithmic Fourier transform (LFT), which they have applied in a geophysical context~\cite{LFT1988}. In its original form however, the LFT is only applicable under very restrictive assumptions on the properties of the function $f(\nu)$ under consideration (e.g. $f(0)=0$) and on the allowed range of the trade-off parameter, which is necessary to adjust the LFT according to the asymptotics of $f(\nu)$.\\
The aim of this work is to present a generalized version of the logarithmic Fourier-Laplace transformation that in particular applies to functions with nonintegrable power-law tails. We give the corresponding definition in section~\ref{sec:def} and show how the original restrictions can be lifted to extend to generalized functions~\cite{GelfandBook}. Moreover, in section~\ref{sec:convergence} we give a proof  that the LFT converges exponentially fast in the number of grid points used for the numerical evaluation, provided the function satisfies certain analyticity conditions. Furthermore, we discuss how the theorem can be applied for practical purposes and in particular show that functions with algebraic tails are perfectly amenable to the LFT. In section~\ref{sec:idealTOP}, we find an ideal set of the trade-off parameters, based on the asymptotic behavior of the input data and extend the excellent performance of the LFT to convolutions in section~\ref{sec:Convo}. In section~\ref{sec:Ex} we provide several classes of mathematical examples highlighting the advantages of LFTs over FFTs and discuss possible optimizations. Finally, we show in section~\ref{sec:physEx} how the LFT can be applied to typical multiscale problems in physics on the example of a density-density correlation function and a simple variant of mode-coupling theory. We conclude in section~\ref{sec:con}.

\section{Definition}\label{sec:def}
\subsection{Mathematical Formulation}
Following the standard convention in the physics literature, we define the Fourier transform of a function $\hat{f}(t)$ in the time domain as
\begin{align}\label{eq:FT}
f(\nu)=\mathcal{F}(\hat{f})(\nu)=\int_{-\infty}^\infty dt \hat{f}(t)e^{i \nu t}\,,
\end{align}
while the inverse transform to frequency $\nu$ is given by
\begin{align}\label{eq:ILFT}
\hat{f}(t)=\mathcal{F}^{-1}(f)(t)=\int_{-\infty}^\infty\frac{d\nu}{2\pi}f(\nu)e^{-i\nu t}\,,
\end{align}
for both $f, \hat{f} \in L^1[\mathbb{R} , \mathbb{C}] $. 
In the following, we utilize the LFT to extend the set of argument functions to  include certain distributions, the precise properties of which we state below. %\red{\emph{We need to be more specific as to which distributions. Maybe we can write all functions with $f(\nu)/f(1/\nu)\to \nu^c$ with $c<0$ as $\nu\to\infty$ and no non-integrable poles on the real axis.}}
We introduce the logarithmic frequency and time coordinates $\omega$ and $\tau$ via
\begin{align}\label{eq:IFT}
\nu=\sigma \bar{\nu} e^{\omega}\quad\text{and}\quad t=\eta \bar{t} e^{\tau}\,,
\end{align}
where $\sigma=\pm 1 = \eta$ are necessary to distinguish between the positive and negative real axis, while the prefactors $\bar{\nu}$ and $\bar{t}$ are required for dimensional purposes and will be set to unity in the remainder of this paper. With these definitions, the inverse Fourier transform~\eqref{eq:ILFT} can be written as a convolution for every $t\in\mathbb{R}$:
\begin{align}\label{eq:convo}
\begin{split}
&\hat{f}(\eta |t|)= e^{-k \tau}\\
& \left. \times\!\!\sum_{\sigma=\pm 1}\! \int\! \frac{d\omega}{2\pi}f(\sigma  e^{\omega})e^{k(\omega+\tau)-i \sigma\eta  \exp{(\omega+\tau)}}e^{(1-k)\omega} \right|_{\tau =\ln |t|}\!,
\end{split}
\end{align}
where $k \in \mathbb{R}$ denotes the trade-off parameter~\cite{LFT1988}.
By the help of the convolution theorem of Fourier analysis (see also Eq.~\eqref{eq:convTheorem} below), the integral in~\eqref{eq:convo} can be reformulated in terms of the product of two Fourier transforms
\begin{widetext}
\begin{align}\label{eq:LFTpre}
\hat{f}(\eta |t|)= \frac{e^{-k \tau}}{2\pi}\sum_{\sigma=\pm 1}\mathcal{F}_{s\to\tau}\left[\mathcal{F}_{\omega\to s}\left(f(\sigma  e^{\omega})e^{(1-k)\omega}\right)(s)\mathcal{F}^{-1}_{x\to s}\left(e^{k x-i\sigma \eta   \exp{(x)}}\right)(s)\right]\left(\tau=\ln|t|\right)\,,
\end{align}
\end{widetext}
provided that $k$ is chosen such that each of the three Fourier integrals converges, the conditions for which we will detail now.\\
Since we ultimately aim for a numerical implementation, the LFT can in general only be applied if 
	\begin{align}
	F_\sigma(\omega):=f(\sigma  e^{\omega})e^{(1-k)\omega} \in L^1\, ,
	\end{align} 
	such that the Fourier transformation 
	\begin{align}
	g_\sigma(s):=\mathcal{F}_{\omega\to s}\left(f(\sigma  e^{\omega})e^{(1-k)\omega}\right)(s)
	\end{align} exists in the integral sense of Eq.~\eqref{eq:FT}.
	Regarding the original function $f(\nu)$ this statement is equivalent to
	\begin{align}\label{eq:integrability}
		\int_{0}^{\infty} d\nu \left|f(\sigma\nu)\right| |\nu|^{-k} < \infty\,.
	\end{align}
 In the particular case of a power-law behavior, i.e. $f(\nu)\to\nu^a$ for $|\nu|\to 0$ and $f(\nu)\to\nu^b$ for $|\nu|\to\infty$, the trade-off parameter has to be chosen according to
\begin{align}\label{eq:condition}
1+b<k<1+a \,.
\end{align}
As a result, for theses functions the LFT even admits a pole of $f$ located at the origin or a branch cut beginning just there, as well as nonintegrable, algebraically growing asymptotics, provided that they can be controlled by an appropriate value of $k$. %\red{In particular, $f(\nu)$ only has to be analytic everywhere along the real axis except at the origin. This exactly corresponds to the typical situation encountered in the context of correlation functions, where poles at positive and negative frequencies are always shifted into the complex plane.}\green{I would keep the linebreak}\\
%\\

Applying the definition of the $\Gamma$ function the $f$-independent inverse Fourier transform in Eq.~\eqref{eq:LFTpre} can be formally rewritten as
\begin{align}
\begin{split}
h_{\sigma\eta}(s):=&\mathcal{F}^{-1}_{x\to s}\left(e^{k x- i\eta \sigma   \exp{(x)}}\right)(s)\\=&\frac{1}{2\pi} \left(i\sigma\eta\right )^{is-k}\Gamma(k-i s)\,,
\end{split}
\end{align}
for $k \in \mathbb{R}\setminus \mathbb{Z}_0^-$, where the exclusion of nonpositive integers is due to the poles of the  Gamma function $\Gamma(k-is)$. 
We point out that this result has to be considered as the analytic continuation of the integral representation
\begin{align}
\int dx\; e^{k x -i \eta \sigma \exp\left(x\right)} e^{- i s x} =\left(i\sigma\eta\right )^{is-k}\Gamma(k-i s) \, 
\end{align}
that, indeed, only holds if $0<k<1$, as emphasized by Haines and Jones~\cite{LFT1988}.\\
Finally, we have to consider the transformation $\mathcal{F}_{s \to \tau}(g_\sigma(s) h_{\sigma\eta}(s))$ from the auxiliary variable $s$ to $\tau$ in Eq.~\eqref{eq:LFTpre}. Since in any practical implementation the factor $g_\sigma(s)$ will only be known in an approximate, discretized form, no analytic continuation can be applied and we have to demand that $g_\sigma \cdot h_{\sigma\eta} \in L^1$.
Given the asymptotics of the product~\cite{frei05book}
	\begin{align}\label{eq:asymptotics}
	\left|\Gamma(k-is)(i\sigma \eta)^{is-k}\right|\!\propto\!
	\begin{cases}
	\!\!\sqrt{2 \pi} |s|^{k-1/2} e^{- \pi |s|}\!\!\!\! & \sigma \eta s \to \infty  \\
	\!\!\sqrt{2 \pi} |s|^{k-1/2}               \!\!\!\!& \sigma \eta s \to -\infty
	\end{cases},
	\end{align}
	we conclude that $g_\sigma \cdot h_{\sigma\eta} \in L^1$ requires $g_\sigma$ to satisfy $\lim_{|s| \to \infty}|s|^{k+1/2} g_\sigma(s)=0$. According to the lemma of Riemann-Lebesgue for differentiable functions~\cite{koer89book}, the latter condition is fulfilled if $F_\sigma(\omega)$ is at least 
\begin{align}
	n:=\max(0, \lceil k+1/2\rceil)
	\label{eq:n}
	\end{align}
	 times differentiable with the derivatives $F^{(\blue{l} )}_\sigma(\omega) \in L_1$, for $0\leq l \leq n$. With respect to the original function $f(\nu)$ this implies that $f^{(n)}(\nu)$ exists, while the integrability condition on $F^{(n)}(\omega)$ reduces to Eq.~\eqref{eq:integrability}, as can be shown by partial integration.

All in all, the logarithmic Fourier transform reads
\begin{widetext}
\begin{align}
\hat{f}(\eta |t|)= \frac{e^{-k \tau}}{(2\pi)^2}\sum_{\sigma=\pm 1}\mathcal{F}_{s\to\tau}\left[\mathcal{F}_{\omega\to s}\left(f(\sigma  e^{\omega})e^{(1-k)\omega}\right)(i\sigma\eta)^{is-k}\Gamma(k-i s)\right]\left(\tau=\ln|t|\right)\,,
\label{eq:LFT}
\end{align}
\end{widetext}
which can be applied with a given value of the trade-off parameter $k \in\mathbb{R} \setminus \mathbb{Z}_0^-$ to all functions $f(\nu)$, that satisfy the summability criterion~\eqref{eq:integrability} and are $n$ times differentiable, with $n$  set by Eq.~\eqref{eq:n}.

The computation of the Fourier transformation $f(\sigma|\nu|)=\mathcal{F}(\hat f(t ))$ via the LFT follows analogously, with the same conditions on $\hat f(t )$. It yields the same result as in Eq.~\eqref{eq:LFT}, yet, with the replacements $f \to \hat{f}$, $\sigma \leftrightarrow \eta$, $\tau \leftrightarrow \omega$ and $(i\sigma\eta)^{is-k}\to (-i\sigma\eta)^{is-k}$, as well as an additional factor of $2\pi $ on the right-hand side of Eq.~\eqref{eq:LFT}.%\\

We remark that the LFT can readily be generalized to Fourier-Laplace transforms of the form
\begin{align}
\mathcal{FL}(f)(t)=\int_{-\infty}^\infty\frac{d\nu}{2\pi}f(\nu)\exp(e^{i\phi}\nu t)\,,
\end{align}
with $\phi \in [0,2\pi[$. The result in~\eqref{eq:LFT} still holds, merely with the factor $(i\sigma\eta)^{is-k}$ substituted by $(e^{i\phi}\sigma\eta)^{is-k}$. Half-sided transforms, which correspond to the standard definition in case of the Laplace transformation, are simply obtained by restricting the outermost sum of Eq.~\eqref{eq:LFT} to $\sigma=1$. Clearly, in most applications Laplace transforms, that is $\phi=\pi$, involve quickly converging integrals. Therefore we will focus on the most critical case of Fourier transforms ($\phi=\pi/2$) and inverse Fourier transforms ($\phi=3\pi/2$), where the transformation kernel entails no exponential suppression of large arguments. Nonetheless, we stress, that even for exponentially decaying integrals the logarithmic transforms are orders of magnitude faster than equidistant grids as is highlighted by the trivial example $f(\nu)=e^{-|\nu|}$ in section~\ref{sec:Ex}.

After having discussed the mathematical framework of the LFT let us briefly comment on the role of $k$ (see also Ref.~\cite{LFT1988}). 
The term trade-off parameter refers to the fact that $k>0$ ($k<0$) suppresses both the integrand of the Fourier transformation $\mathcal{F}_{\omega \to s}$ in the definition~\eqref{eq:LFT} (cf. also Eq.~\eqref{eq:integrability}) for large $\omega \to \infty$ (small $\omega \to -\infty$) and the result in the $\tau \to \infty$ ($\tau \to -\infty$) limit, which corresponds to $t \to \pm \infty$ ($t \to 0$), due to the overall prefactor. Simultaneously, the convergence in the opposite limits is diminished. This dependence on $k$ can be utilized to tune the properties of the LFT to suit the asymptotics of $f$.\\
A different perspective on the LFT is opened by the interpretation of the trade-off parameter as a shift of the final integration over $s$ to a contour in the complex plane. The discussion of the LFT in terms of contour integrals, which are sensitive to the analytic structure of $f$, is crucial to understand the convergence properties of the LFT detailed in section~\ref{sec:convergence}

\subsection{Numerical implementation}
So far, we have only utilized exact analytical reformulations of the problem. However, upon introducing exponential grids with index $\;n\in\{1,2,...,N\}$ in both frequency and time
\begin{align}\label{eq:grids}
\begin{split}
\nu_n&= e^{\omega_n} \quad \text{and} \quad \nu_{-n}=- e^{\omega_n} \quad\text{with} \quad\omega_n=\Delta \omega \left(n +\omega_s\right)\;\\
t_n&= e^{\tau_n} \quad \,\text{and} \,\quad t_{-n}=- e^{\tau_n}\,\quad \text{with}\,\quad \tau_n=\Delta\tau \left(n +\tau_s\right)
\end{split}
\end{align}
and discretizing the auxiliary space via $s_n=\Delta s( n + s_s)$, the usefulness of the form~\eqref{eq:LFT} is immediately revealed: The equidistant grids in $\tau$ and $\omega$ make it possible to take advantage of the efficient FFT algorithm -- even for Laplace transforms, where fast algorithms otherwise require more elaborate methods from approximation theory~\cite{Rokhlin1988, Strain1992} -- while covering low frequencies and short times with a high density of points, as opposed to a reduced sampling density at large arguments. Since in many physical applications the high-energy or frequency range shows an algebraic behavior, this covering of the frequency (momentum) and time (position) domain will be very favorable under many circumstances.
Important physical examples include generic correlation functions in frequency and momentum space, while in a critical theory algebraic tails appear in the position and time argument. For instance the momentum distribution $n(k)$ of ultracold Fermions in the vicinity of an open-channel dominated Feshbach resonance~\cite{chin10Feshbach} obeys the Tan energy theorem~\cite{Tan08energy}: For momenta $k$ that exceed any intrinsic inverse length scale $n(k)$ decays like $\mathcal C/k^4$, where $\mathcal C$ is the observable Tan contact density~\cite{kuhn11contact, hoin13contact}. On the other hand, the phase transition to the superfluid is signaled by an instability of the pair propagator in the low-momentum limit. 
In this system, the LFT has been applied to study the phase diagram in the presence of a finite spin imbalance~\cite{Frank2018}.
Furthermore, similar challenges arise in the efficient simulation of analog low/high pass filters~\cite{Christensen1990}, in the context of signal processing as well as in the numerical solution of differential equations~\cite{Boyd2001}.\\
In addition to the convenient distribution of points, the grid~\eqref{eq:grids} acquires a high degree of flexibility as the step sizes $\Delta \omega, \Delta \tau$ and $\Delta s$ and the linear shifts $\omega_s$ and $\tau_s$, that play the role of the prefactors $\bar{t}$ and $\bar{\nu}$, together with $s_s$ can be chosen at will. This allows to adjust the method to the asymptotics of various functions as is shown in section~\ref{sec:Ex}.
The standard choice for all the shifts is $-N/2$ in order to cover positive and negative exponents equally. We will return to the question of how to determine the ideal transformation parameters in section~\ref{sec:idealTOP}. 

Before continuing, we remark however, that functions with important features on intermediate scales which cannot be considered as part of the asymptotics of small or large arguments (not even by using the entire set of parameters available in Def.~\eqref{eq:grids}), will yield no advantage over an ordinary FFT. Such cases appear for double-peak structures whose centers are too far apart to be scaled to the high grid-density at $\omega \to 0$ without including an impractically large $N$. Similarly, functions that oscillate uniformly on all scales with a fixed frequency $\bar{\omega}$ will be inevitably undersampled by the given grid at frequencies $\omega \gtrsim \ln(\bar{\omega}/\Delta\omega)$.

\section{Convergence properties}\label{sec:convergence}
\subsection{Theoretical perspective}\label{sec:proof}
Now we address the issue of how efficiently a function $f(\nu)$ that obeys the properties stated below Eq.~\eqref{eq:LFT} can be sampled and then Fourier transformed on the exponential grid. This requires to answer the question of how quickly the sum
\begin{align}\label{eq:LFTsum}
\begin{split}
\hat{f}_N(t_{\eta n})=\; & e^{-k \tau_n}\!\!\sum_{\sigma=\pm 1}\! \sum_{l=1}^N\frac{\Delta s}{2\pi}e^{i s_l \tau_n}(i\sigma \eta  )^{i s_l -k}\Gamma(k-is_l)\\
&\! \!\qquad \qquad\cdot\sum_{m=1}^N \frac{\Delta \omega}{2\pi} f(\sigma  e^{\omega_m})e^{(1-k)\omega_m}e^{i\omega_m s_l} \,
\end{split}
\end{align}
representing the numerical, discrete approximation on the grids defined in Eq.~\eqref{eq:grids}
converges towards the exact integral~\eqref{eq:LFT} as $N \to \infty$. First of all, we note that Eq.~\eqref{eq:LFTsum} indeed approaches the LFT from Eq.~\eqref{eq:LFT}. To see this one has to consider the limits of the largest values $|\omega_{\pm N}|, |s_{\pm N}| \to \infty$ at vanishing stepsizes $\Delta \omega, \Delta s \to 0.$ Taking the latter limit yields well-defined integrals on finite intervals, since all terms represent measurable functions. In particular, the sum in the second line can be interpreted as Fourier coefficient $F_l$ of the periodic function $F_\sigma(\omega)$ with period $2N\Delta\omega$. The differentiability of $F_\sigma(\omega)$ then implies the asymptotic behavior $F_l \lesssim C s^{(n+1)}$ with a positive constant $C$~\cite{koer89book}, such that the limit $|\omega_{\pm N}|, |s_{\pm N}| \to \infty$ exists and by its uniqueness we recover the definition of the LFT.

%\red{Returning to the question how fast the convergence can be made, we proof below that we obtain exponential convergence in the number of grid points, provided the function $f(\sigma e^\omega)$ satisfies certain analyticity conditions. In the next subsection we show that many of the functions relevant in physical applications fulfill these and can therefore be treated very efficiently with the LFT at an extraordinarily high precision.}\\
Beyond the mere existence, we now show that under conditions satisfied in many relevant application, the LFT converges exponentially fast in the number of grid points.\\

	\textbf{Theorem:} Let $f(\sigma e^\omega)$ be a function that is analytic in a closed strip of width $R^{(1)}>0$ around $\bar{\mathbb{R}}$, i.e. the affinely extended real axis of the logarithmic argument $\omega$, and whose asymptotic behavior can be controlled by a suitable choice of the trade-off parameter $k$, such that $F_\sigma(\omega)\in\mathcal{S}(\mathbb{R})$, where $\mathcal{S}$ denotes the space of Schwartz functions~\cite{GelfandBook}. Then the deviation of the approximation Eq.~\eqref{eq:LFTsum} from the exact expression Eq.~\eqref{eq:LFT} vanishes exponentially in the number of grid points $N$. \\
	
	\textbf{Proof:}
The rate of this convergence will not depend on the exact values of the centers of the grids $\omega_s$, $\tau_s$ and $s_s$. To keep the notation simple, we will in the following assume them to be given by integers. 

\begin{figure}[htp]
\begin{center}
\includegraphics[width=0.9\columnwidth]{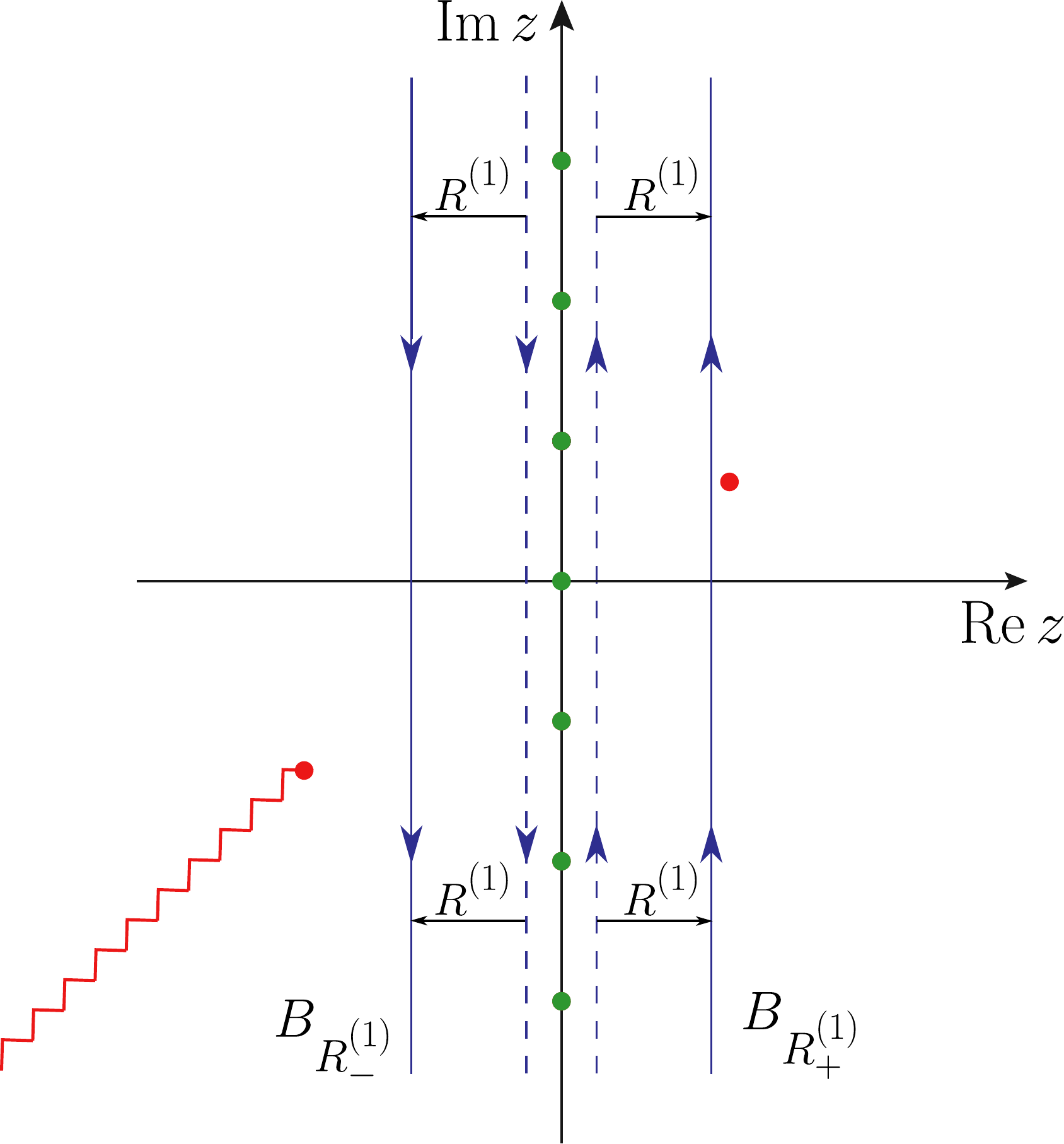}
\caption{(Color online) The original sum in $S_1^\sigma(s)$, evaluated at the green dots, is replaced by an integration contour along the imaginary axis (blue dashed line) and then shifted by a finite real part $\pm R^{(1)}$. The closest approach of a non-analyticity of $f(\sigma e^{-iz})$ -- here symbolized by a red dot for a pole and a red zigzag line for a branch cut -- to the imaginary axis determines $R^{(1)}$.}
\label{fig:Boundary}
\end{center}
\end{figure}

Obviously, $ \mathcal{S} \subset L^1$ and the Schwartz functions satisfy the differentiability condition~\eqref{eq:n} by definition. 
Furthermore, the integral  
\begin{align} \label{eq:I1}
I^\sigma_1(s) = \int \frac{d\omega}{2\pi}f(\sigma e^{\omega})e^{(1-k)\omega}e^{i s\omega} \, ,
\end{align}
is finite and itself a Schwartz function since the integrand is an element of $\mathcal{S}(\mathbb{R})$. By virtue of the Paley-Wiener theorem~\cite{titc75book} $I^\sigma_1(s)$ is analytic in a strip around the real $s$-axis, whose width $R^{(2)}>0$ is determined by the asymptotic decrease of $F_\sigma(\omega)$, which at least is exponential.
%https://math.stackexchange.com/questions/138302/compactly-supported-function-whose-fourier-transform-decays-exponentially
Furthermore, the truncation error due to the finite summation interval, scales like $|F_\sigma(\omega_{\pm N})|$ and thus in any case merely gives rise to exponential corrections.
Therefore, we consider right away the infinite sum. The latter can be replaced by a contour integral around the imaginary axis in the mathematically positive direction, which reads
\footnote{This method of rewriting sums in terms of contour integrals is also used to evaluate Matsubara sums in the context of finite temperature quantum field theory~\cite{abri75, altl2010book,fett71}}
\begin{align}\label{eq:S1}
\begin{split}
& S_1^\sigma(s) =\sum_{m \in \mathbb{Z}} \frac{\Delta\omega}{2\pi}f(\sigma e^{\omega_n})e^{(1-k)\omega_n}e^{i s\omega_n} \\ & =
\begin{cases}
\oint\frac{d z}{2\pi i}f(\sigma e^{-i z})(1+n_B(z))e^{-(1-k+is)iz} &  \Re(s)<0 
\\
\oint\frac{d z}{2\pi i}f(\sigma e^{-i z})(n_B(z))e^{-(1-k+is)iz}   &  \Re(s)>0
\end{cases} \,.
\end{split}
\end{align}
Here $n_B(z)=1/(\exp(2\pi z/\Delta\omega)-1)$ is the Bose-Einstein distribution with "inverse temperature" $\beta_\omega=2\pi/\Delta\omega$, whose simple poles at $\beta_\omega \omega_n$ make sure that one recovers the original series with the help of the residue theorem. 
Subtracting the exact integral $I_1^\sigma(s)$, which is
also taken along the imaginary axis, one obtains for the difference

\begin{align} \label{eq:E1}
\begin{split}
&E_1^\sigma(s)= S^\sigma_1(s)-I^\sigma_1(s)=\\& =\; \int_{R^{(1)}-i \infty}^{R^{(1)}-+ i \infty} \frac{d z}{2\pi i}f(\sigma e^{-i z})(n_B(z))e^{-(1-k+is)iz} 
\\ &+  \;\int_{-R^{(1)}+i \infty}^{-R^{(1)}- i \infty} \frac{d z}{2\pi i}f(\sigma e^{-i z})(1+n_B(z))e^{-(1-k+is)iz} \, ,
\end{split}
\end{align}
where we have made use of Cauchy's theorem to deform the integration contour such that it remains within boundary of the analytic domain of $f(\sigma \exp(\omega))$ (cf. Fig.~\ref{fig:Boundary}).
Extracting the dominant exponential behavior we can write for the error
\begin{align}\label{eq:E1final}
E_1^\sigma(s)= e^{R^{(1)}(-\beta_\omega+|s|)} F^\sigma_1(s,R^{(1)})\, .
\end{align}
The function $F_1$ arises from the remaining integrals and is the Fourier transformation of an integrable function, which in particular implies that the exponential prefactor indeed yields the leading asymptotic behavior for $|s| \to \infty$ due to the lemma of Riemann and Lebesgue. Moreover, by increasing $\beta_\omega \sim N$ the error becomes exponentially small in the number of grid points, provided $|s|<\beta_\omega$. From a theoretical perspective, we can take the limit $N \to \infty$ and achieve exponential convergence uniformly in $s$.

Next, one has to investigate the convergence properties of the remaining sum over the auxiliary variable $s$ in~\eqref{eq:LFTsum}. 
First, we focus on the properties of the exact integral
\begin{align}\label{eq:I2}
I_2^\sigma(\tau) = \int \frac{ds}{2 \pi} \Gamma(k-is)(i\sigma \eta)^{is-k} I_1^\sigma(s) e^{i s\tau} \,  .  
\end{align}
The asymptotics of the product~\eqref{eq:asymptotics}
	%~\cite{frei05book}
%\begin{align}\label{eq:asymptotics}
%\left|\Gamma(k-is)(i\sigma \eta)^{is-k}\right|\propto
%\begin{cases}
%\sqrt{2 \pi} |s|^{k-1/2} e^{- \pi |s|}\!\!\!\! & , \sigma \eta s \to \infty  \\
%\sqrt{2 \pi} |s|^{k-1/2}               \!\!\!\!& , \sigma \eta s \to -\infty
%\end{cases}
%\end{align}
gives at most rise to an algebraic growth, while for $|s| \to \infty$ the integral $I_1^\sigma(s)\in\mathcal{S}$ decays exponentially fast, 
thus rendering $I_2^\sigma(\tau)$ well-defined.
To estimate the error arising from replacing the analytic expression in Eq.~\eqref{eq:LFT} by a discretized numerical approximation we consider the difference
\begin{align}\label{eq:E2}
\begin{split}
E^\sigma_2(\tau)=&I^\sigma_2(\tau)\\&-\sum_{|l|\leq N} \frac{\Delta s}{2\pi}\Gamma(k-is_l)e^{i s_l \tau}(-i\sigma \eta)^{i s_l-k}S^\sigma_1(s_l) \,  .
\end{split}
\end{align}
Since the difference $E_1^\sigma(s) = S_1^\sigma(s) -I_1^\sigma(s)$ becomes exponentially small with decreasing $\Delta\omega$, we can replace the sum by the exact integral $I_1^\sigma(s)$. Furthermore, the sum can be extended to $l \in \mathbb{Z}$ because the truncation to $|l| \leq N$ neglects only terms that are exponentially suppressed due to the asymptotics of $I_1^\sigma(s)$ that, provided $s_{\pm N}$ are large enough, determines the exponential tails of $S_1^\sigma(s)$. Then, the error $E_2^\sigma(\tau)$ can be treated analogously to Eq.~\eqref{eq:E1} in terms of complex contour integrals by introducing the Bose distribution $n_B(z)$, which now involves the inverse temperature $\beta_s=2\pi/\Delta s$. In this step we have to analytically continue the integral $I_1^\sigma(s \to z)$, which is in general not known in closed form for genuine complex arguments $z$. Yet, performing this continuation numerically is not required as we need the expression only on a formal level to determine the error.
Like above, we shift the contour integrals into the complex plane to the fixed finite real parts $\pm R^{(2)}$ within the width of the analytic domain of $I_1^\sigma(z)$, but have to take into account the simple poles of the $\Gamma$ function located at the nonpositive integers enclosed by the modified contour. Altogether, we can estimate the error by
\begin{align}\label{eq:E2new}
E^\sigma_2(\tau)= e^{R^{(2)}(-\beta_s +|\tau|)}F_2^\sigma(\tau, R^{(2)})+E^\sigma_{2,\Gamma}(\tau)\, .
\end{align}	
Once again $F_2^\sigma(\tau, R^{(2)})$ is the Fourier transform of an integrable function and does not overcome the leading exponential such that the first term vanishes exponentially for all $|s|$ in the limit $\beta_s \sim 1/\Delta s \sim N\to\infty$. The second term summarizes the contributions from the residues $\text{Res}(\Gamma,-m) = (-1)^m/m!$ for $m \in \mathbb{N}_0$ and reads
\begin{align}\label{eq:E2Gamma}
\begin{split}
&E^\sigma_{2 \Gamma}(\tau) =  -\!\!\sum_{\substack{m=\lceil -k \rceil\\ \land m \geq 0}}^{\lfloor R^{(2)}-k\rfloor} \frac{(-1)^m}{m!} \frac{e^{(k+m)\tau}(i \sigma \eta)^m I_1^\sigma(-i(m+k))}{e^{\beta_s(k+m)}-1} \\
&\quad\;\; -\!\! \sum_{\substack{m=\lceil -R^{(2)}-k \rceil\\ \land m \geq 0}}^{\lfloor -k \rfloor} \frac{(-1)^m}{m!} \frac{e^{(k+m)\tau}(i\sigma \eta)^m I_1^\sigma(-i(k+m))}{1-e^{-\beta_s(k+m)}}\, .
\end{split}
\end{align}
Note that the Bose functions control the exponential function $\exp((k+m)\tau)$ via $\beta_s$ in analogy to  the first term in Eq.~\eqref{eq:E2new}. 

In total, we observe that both $E_1^\sigma$ and $E_2^\sigma$ decrease exponentially in the limit $\Delta_\omega,\Delta s \to 0$ for all $s$, while all intermediate steps do not violate this scaling. $\square$\\
Note that for $\phi\in]\pi/2,3\pi/2[$, which includes the half-sided Laplace transform, convergence is even faster. This, however, takes its toll, when considering the inverse transform, where the exponential growth of $|\Gamma(k-is)|^{-1}\sim e^{\pi|s|/2}$ in the case of the Laplace transform severely magnifies errors and limits the useful interval in $s$ and thereby the attainable precision. For the special case of $k=1/2$, this has been analyzed in detail by Epstein and Schotland~\cite{Epstein2008}, who, given a noise $\delta$ on the input, also derive a bound on the maximum resolution $\Delta\omega$ and precision $\epsilon$ of the numeric inverse Laplace transform:
\begin{align}\label{eq:laplace}
1\leq\frac{2\Delta\omega}{\pi^2}\ln{\left(\sqrt{2\pi}\frac{\epsilon}{\delta}\right)}\;.
\end{align}
As we will see below, this bound can be significantly improved by the use of the theorem in Sec.~\ref{sec:proof} in combination with the knowledge of $R^{(1)}$.
%After subtraction of the corrections due to $E_{2\Gamma}^\sigma$, the error is a flat function in s-space with amplitude $\simeq \epsilon$. It has to be compared with the initially transformed function $S_1^\sigma(s)\propto e^{-R^{(1)}(\beta_\oemga-|s|)}+e^{-R^{(1)}|s|}$. The minimal ratio determines the error in $\nu$ (up to the slope set by $k$).
\subsection{Practical aspects}
Regarding the implementation of the LFT for practical purposes, it is helpful to relate the analytic properties of the function $f$ in the logarithmic argument $\omega$ to the original argument $\nu$, since $f(\nu)$ corresponds to the form that is given in most applications. From this we find bounds for $R^{(1)}$ for some relevant classes of functions. First of all, the analyticity of $f(\sigma e^\omega)$ implies that $f(\nu)$ is also analytic for $\nu \in \mathbb{R}\setminus\{0\}$.
Moreover, the exponential decay of $F_\sigma(\omega\to\pm\infty)$ requires that a constant $c_\sigma>0$ exists, such that
\begin{align}\label{eq:ana_asympt}
\left|f(\sigma|\nu|)f\left(\frac{\sigma}{|\nu|}\right)\right|\underset{\nu\to 0}{<}|\nu|^{c_\sigma}\;.
\end{align}
To connect the analyticity properties of $f$ with respect to the variables $\nu$ and $\omega$, we first note that the analytic strip of width $R^{(1)}$ in $\omega$ is equivalent to the statement that for each $\omega_0 \in \bar{\mathbb{R}}$ there is an $R^\sigma_{\omega_0}>0$, such that the Taylor series of $f(\sigma e^\omega)$
\begin{align}\label{eq:Taylor}
	f(\sigma e^\omega)= \sum_{n=0}^{\infty} \left(\frac{d^n}{d \omega^n} f(\sigma e^\omega)\right)_{\omega_0} \frac{(\omega-\omega_0)^n}{n!}
\end{align}
 converges for all $|\omega-\omega_0| \leq R^\sigma_{\omega_0}$. The width of the strip is given by 
\begin{align}\label{eq:R1}
R^{(1)} = \inf_{\substack{\omega_0 \in \bar{\mathbb{R}}\\ \sigma= \pm 1}} R_{\omega_0}^\sigma\, .  
\end{align}
To estimate the minimal size of the analytic domain of $f(\nu)$, we first determine the image of the line $\omega_0 + i \lambda R^{(1)} $, with $\lambda \in [-1,1]$ centered around the real $\omega_0$ under the mapping~\eqref{eq:IFT}. This yields circular sectors of radius $\nu_0=\exp{\omega_0}$ that are symmetric around the real $\sigma \nu$ half-axis and centered at $\nu=0$. The half opening angle is given by $\min(\pi/2,R^{(1)})$. The restriction of the angle arises from the separation of the variable $\nu$ with respect to $\sigma = \text{sign}( \Re(\nu))$. Taking the union over all $\omega_0 \in \mathbb{R}$, which yields the domain of analyticity of $f(\nu)$ gives rise to an infinite cone $|\Im(\nu)|<\tan{(\theta)}|\Re(\nu)|$ with opening angle
	\begin{align}
		\theta(R^{(1)}) = \min(\pi/2,R^{(1)}) >0 \, .
	\end{align} 
Note that for asymptotically large linear-frequency arguments the required width in $\nu$ space grows linearly, but admits nonanalyticities close to the origin. This is to be expected since nonsmooth variations of $f(\nu)$ that are related to nonanalyticities are very well captured by the exponentially dense grid in the vicinity of the origin. On the other hand, the same behavior at large frequencies will cause deteriorated numerical results due to severe undersampling of sharp features.\\
In view of these arguments it becomes apparent that 
 algebraic functions with the asymptotic behavior $f(\nu \to 0)\sim \nu^a$ and $f(\nu \to \infty)\sim \nu^b$ with $a>b$ are ideal candidates for the application of the LFT, since they trivially satisfy condition~\eqref{eq:ana_asympt} and the transformation $\nu \to \sigma \exp(\omega)$ removes any nonanalyticity located at the origin, which arises from any $a \in \mathbb{R}\setminus \mathbb{N}_0$. We recall that the above considerations do not make any reference to the integrability properties of $f(\nu)$, which in case of algebraic functions can be controlled by the trade-off parameter. 
 
 Another class of asymptotic behavior is given by exponential functions, which we discuss with the help of the simple example of a single, dominant exponent $f(\nu) \sim \exp(\alpha (\sigma \nu)^c)$ for $\sigma \nu \to 0, \infty$, with $c \in \mathbb{R} \setminus \{0\}$ and $\alpha \in \mathbb{C}\setminus\{0\}$ to include oscillating functions. The prefactor is allowed to contain an arbitrary algebraic function to which the problem would be reduced for $\alpha=0$ or $c=0$. 
 After the mapping to $\omega = \Re \omega + i \Im \omega$ we have
 \begin{align}
 	\left|e^{\alpha(\sigma \nu )^c}\right| = \left| e^{\exp(c\, \Re \omega )\left(\Re \alpha \cos(c\, \Im \omega )- \Im \alpha \sin(c\,\Im \omega)\right) }\right| \, .
 \end{align}  
 For real $\omega$ the limit $c\, \Re \omega \to \infty$ requires $\Re \alpha <0$, since otherwise the inner Fourier transform $\mathcal F_{\omega \to s}$ of the LFT in Eq.~\eqref{eq:LFT} is not defined. Note that the trade-off parameter cannot be used  to remedy the super-exponential growth for $\Re \alpha >0$, which is also beyond the scope of the notion of generalized Fourier transformations~\cite{GelfandBook}. Furthermore, we observe that even if $\Re \alpha <0$, the boundary of the analytic strip cannot overcome the constraint
 \begin{align}
 	R^{(1)}  < \left| \frac{1}{c} \arctan \left(\frac{\Re \alpha}{\Im \alpha}\right) \right| 
 \end{align}
 because $f(\sigma e^\omega)$ diverges exponentially for $\Im \omega$ beyond that value if $c\, \Re \omega \to \infty$. Considering an oscillatory example $f(\nu) \sim e^{i \nu}$, we obtain $R^{(1)}=0$, irrespective of an integrable algebraic prefactor. Therefore, the convergence of the LFT is degraded to an algebraic one on a fundamental level as mentioned at the end of Sec.~\ref{sec:def}.\\
%\begin{align}
%\begin{split}
%\sigma \nu &= \sigma \nu_0 + R^\sigma_{\nu_0} e^{i \varphi}\, , \text{if } |\nu_0| \to \infty \\
%(\sigma \nu)^{-1} &= (\sigma \nu_0)^{-1} + R^\sigma_{\nu_0} e^{i \varphi} \, , \text{if } |\nu_0| \to 0 \, ,
%\end{split}
%\end{align} 
%where $\varphi \in [0, 2\pi)$, we conclude for the extent $R^\sigma_{\nu_0}$ of the analytic region of $f(\nu)$ around $\sigma \nu_0$  
%\begin{align}
%	R^\sigma_{\nu_0} = \left(e^{R_{\omega_0}}-1\right) \begin{cases}
%	\sigma \nu_0\, , & \text{if } \sigma \nu_0 \to \infty \\
%	(\sigma \nu_0)^{-1}\,, & \text{if } \nu_0 \to 0 \, .
%	\end{cases}
%\end{align} 
More generally the width of the analytic strip $R^{(1)}$ can be determined from the positions of singularities of $f(\nu)$ in the complex plane.
Modeling $f(\nu)$ in the vicinity of a nonanalyticity at $\hat \nu = \nu_0 \pm i R^\sigma_{\nu_0}$, where $\nu_0 \in \mathbb{R} \setminus \{0\}$, $\sigma = \text{sgn}(\nu_0)$ and $R^\sigma_{\nu_0} >0$, in the form
	\begin{align}
	f(\nu) = \lambda (\nu-\nu_0 \mp i R^\sigma_{\nu_0})^\alpha \, ,
	\end{align} 
	with $\lambda \in \mathbb{C}$ and $\alpha \in \mathbb{R} \setminus \mathbb{N}_0$, we find for the $k^{\text{th}}$ derivative 
	\begin{align}
		|f^{(k)}(\nu_0)| = |\lambda| \prod_{j=0}^{k-1}(\alpha-j) |R^\sigma_{\nu_0}|^{\alpha-k} \, .
	\end{align}
This form can be used to extract $R_{\nu_0}^\sigma$ and therefore $R^{(1)}$ and is accessible even for numerical data via finite difference approximations.\\
Regarding the application of the LFT in practice, we note that it can be implemented with any existing library of FFTs. However, one additionally has to evaluate the Gamma function for complex arguments. Fortunately, these do not depend on $f$ and can be tabulated if several transformations have to be performed. In any case, all modern programming languages include fast algorithms (e.g. Spouge's approximation~\cite{spou94}) to compute $\Gamma$ and the LFT will, therefore, not be severely slowed down compared to an FFT with the same number of data points.
 
\section{Optimal parameter choices}\label{sec:idealTOP}
So far, we have shown that the LFT converges exponentially fast towards the exact Fourier transform if the analytic structure of $f$ satisfies the conditions of the above theorem. However, we have not yet made any rigorous statements regarding how many points have to be used to reach the desired precision or how the trade-off parameters should be adjusted to obtain an optimal convergence. As we have seen in section~\ref{sec:proof}, the differences between the exact Fourier transform and the LFT approximant are well known, such that generic statements about the ideal parameter choices, as well as significant improvements to the results, can be made.\\
In the following, we will focus mainly on algebraic functions, whose convergence properties can be influenced by the trade-off  parameter in contrast to exponential functions.
First, we recall that the only restriction on $k$ is given by the convergence of the inner integral in the LFT~\eqref{eq:LFT}, thus the trade-off parameter has to be chosen according to $1+b<k<1+a$ and $\notin \mathbb Z_0^-$. From a numerical perspective one has to keep in mind that $k$ has to avoid nonpositive integers by a finite margin to bypass the poles of the $\Gamma$ function. In practice, a deviation of 0.01 turns out to be sufficient. Apart from this constraint it is desirable to choose $k$ close to $k_{\text{opt}}=1+(a+b)/2$ that symmetrizes the asymptotic behavior of the integrand in $I^\sigma_1(s)$ (see Eq.~\eqref{eq:I1}) on both ends of the $\omega$-interval, such that the smallest truncation errors are achieved for the standard value for the centers of the grids (see Eq.~\eqref{eq:grids}). With this trade-off parameter, using $N$ data points, a truncation error no larger than $\epsilon$ requires
\begin{align}\label{eq:cond1}
\Delta\omega=\frac{4}{(b-a)N}\ln(\epsilon)\,.
\end{align}
%In $s$ space the expression
%\begin{align}
%\left|\Gamma(k-is)(i)^{is-k}\right|+\left|\Gamma(k-is)(-i)^{is-k}\right|\sim e^{\left(\frac{1}{2}-k\right)\ln(s)}
%\end{align}
%is, to logarithmic accuracy, of order unity, confer also equations \eqref{eq:asymptotics}.
As argued below equation~\eqref{eq:asymptotics} the decay of the integrand in~\eqref{eq:I2} is dominated by $S_1^\sigma(s)$, the asymptotic behavior of which can be estimated from the contour integral~\eqref{eq:S1}, which gives rise to the asymptotic behavior $S_1^\sigma(s)\sim e^{-R^{(1)}|s|}$. However, approximately at $s_\text{max}=\pi/\Delta\omega$ this function drops below the error $E_1^\sigma(s)$ from Eq.~\eqref{eq:E1final}. 
%As these are well known, they can be subtracted one after the other thereby expanding the accessible interval in $s$ to at least $]-3/2s_\text{max},3/2s_\text{max}[$. Since the gains due to this subtraction are fairly small, we will not rely on it, and instead keep the procedure as simple as possible. 
To compute the remaining Fourier transform with truncation errors that are consistent with the previous steps one consequently demands, that $S_1^\sigma(s_\text{max})<\epsilon$ from which one concludes
\begin{align}
\Delta\omega=-\frac{\pi R^{(1)}}{\ln(\epsilon)}\,.
\end{align}
Together with~\eqref{eq:cond1} this fixes the lower bound of points necessary for an absolute accuracy of roughly $\epsilon e^{-k_\text{opt}\tau}=\epsilon e^{-(1+(a+b)/2)\tau}$ to
\begin{align}\label{eq:minN}
N=\frac{4}{\left(a-b\right)\pi R^{(1)}}\ln^2(\epsilon)\,.
\end{align} 
Remarkably this scales only \emph{logarithmically} with the desired precision. We emphasize that this statement holds in general, even if the optimal choice of the trade-off parameter is prohibited. In this case only the prefactor increases. However, one might suspect that $N$ will be drastically increased by the requirement $\Delta s \ll 1$ in order to control the error $E_2^\sigma(\tau)$, according to equation~\eqref{eq:E2new}. Since the closest non-analyticity of $S_1^\sigma(z)$ to the real axis appears at a distance $R^{(2)} \simeq\min\{k-1-b,1+a-k\}$, we infer from the asymptotics $I_2^\sigma(\tau) \sim \exp(-R^{(2)}\tau)$ and $E_2^\sigma(\tau)\sim \exp(-(\beta_s + |\tau|)R^{(2)})$ that it makes only sense to include values $
|\tau_n| \leq \tau_{\max} = \pi/\Delta s =\beta_s/2$.Therefore, $\Delta s \approx (b-a)\pi/\ln{(\epsilon)}$ is possible. \\
In addition, the error $E_{2\Gamma}^{\sigma}(\tau)$ due to the proximity between the integration contour and poles of the $\Gamma$ function does not influence $\Delta s$ and apart from a prefactor their $\tau$-dependence is exactly known, as can be seen from Eq.~\eqref{eq:E2Gamma}. Thus, if the numerical error in the time domain is dominated by these contributions one can fit the residues in $E^\sigma_{2 \Gamma}$ to the limits $\tau \to \tau_{\pm N}$. In these regimes only numerical noise remains, because by virtue of the lemma by Riemann and Lebesgue the exact function $I_2^\sigma(\tau)$ has decreased below the desired precision threshold. This procedure works particularly well, since the exponential terms are known exactly (see also the examples in section~\ref{sec:Ex}). Moreover, depending on $k$, only a few dominant terms have to be subtracted, while the remaining terms are negligible due to their strictly monotonically decreasing exponents. 
%can be neglected. The reason for doing so is that apart from a simple prefactor the shape of these contributions are exactly known. Thus, for most functions one can easily fit the first few summands of the first line in $E_2(\tau)$ without any loss in precision and subtract them to further improve the quality of the result (see section \ref{sec:Ex}).   

Note that this discussion does not include round-off errors. These will give rise to an additional limitation of the attainable precision, since the finite accuracy of the internal numerical operations sets a bound to the possible precision of the LFT and in particular determines how well the decay of $I_2^\sigma(\tau)$ for $|\tau| \to \infty$ can be resolved. Furthermore, the final multiplication with $e^{-k\tau}$ also affects the error estimate. For negative values of $k$ exceptional precision can be achieved in the regime of small $\tau$. For strongly negative values of $k$, however, these come at the price of enhanced errors at large arguments. 
As long as round-off errors are ignored, these are a minor problem, since the final Fourier transform in~\eqref{eq:IFT} will typically decay to zero at large arguments, which allows to remove the previously discussed systematic errors from $E_{2\Gamma}^\sigma$ (cf. section~\ref{sec:Ex}). In practice, round-off errors, unfortunately, dominate $I^\sigma_2(\tau)$ at large arguments, which sets a lower boundary to the useful interval of the trade-off parameter. Positive values of $k$ on the other hand, which are necessary to treat non-integrable functions, that is those with $b>-1$, will result in undesirably enhanced errors near the origin in the image space (here $t$). Removal of these errors will typically involve fitting the asymptotic behavior to the numerical data within the range of $t$-values, where $\hat{f}(t)>\epsilon e^{-k\tau}$ can be satisfied and extrapolating it towards $t \to 0$. The dynamical compression, that is the diminishing length of this $t$-interval as $k$ increases, is the price to pay for numerically Fourier transforming non-integrable functions.\\
If more than just the minimal number of data points necessary for a given precision are available, one can use them to compensate the enhanced truncation errors and perform several transformations with different trade-off parameters. Larger values of $k$ increase precision for $\tau>0$, while smaller values enhance the accuracy at negative $\tau$. If $k$ can be varied over a wide interval without creating too large truncation errors, this procedure can, for example, be used to mitigate the impact of dynamical compression (see section~\ref{sec:Ex}). To improve the results further one can use $\omega_s$ to shift the list of $\omega_n$ points to optimally sample the asymptotics of the function. For an arbitrary value of $k$ that is compatible with the convergence requirements the condition $F_\sigma(\omega_{\pm N})<\epsilon$ translates to 
	\begin{align}
	\begin{split}
	\omega_s & = \frac{\ln(\epsilon)}{\left(a+1-k\right)\Delta \omega} \\
	N & = -\frac{a-b}{\left(b+1-k\right)\left(a+1-k\right) \pi R^{(1)}}\ln^2(\epsilon)\, ,
	\end{split}	
	\end{align}
which reduces to equation~\eqref{eq:minN} and $\omega_s=-N/2$ if $k=k_\text{opt}$.
Similarly, $s_s$ can be used to find the best distribution of the auxiliary space points $s_n$.\\
Finally, if a precision close to the round-off limit is required, $e^{-k \tau}$ cannot be orders of magnitude larger than $\hat{f}(\eta e^\tau)$ but instead should stay as close as possible to $-\ln(\hat{f}(\eta e^\tau))/\tau$ in the range of $\tau$ arguments of interest.\\
In case of the inverse Laplace transform $\hat{f}(t)\to f(\nu)$ we can use the same procedure for the optimization of the transformation parameters on an input with multiplicative noise of amplitude $\delta$. Given a good estimate of $R^{(1)}$, which equals the width of the analytic strip of $\hat{f}(e^{\tau})$ minus $\pi/2$, this allows us to enhance the relative precision of the result near $\nu=1$ to $\delta^{R^{(1)}/(\pi/2+R^{(1)})}$, where, as opposed to Eq.~\eqref{eq:laplace}, no constraint on $\Delta\omega$ is required. Note that only for $R^{(1)}\ll 1$ one is struck by the dreaded exponential enhancement of errors.

\section{LFT Convolutions} \label{sec:Convo}
One of the most important applications of Fourier transforms in theoretical physics relies on the efficient calculation of convolutions. Due to the convolution theorem
\begin{align}\label{eq:convTheorem}
\int\frac{d\nu'}{2\pi}f(\nu')g(\nu-\nu')=\mathcal{F}(\hat{f}(t)\hat{g}(t))(\nu)
\end{align}
the computational complexity $\propto N^2$ of the direct discretized evaluation of the integral can be reduced to $\propto N \log_2(N)$ when using the FFT algorithm. There are, however, many situations, where the convolution theorem may not be utilized. For example, if one of the factors $f( \nu)$ or $g( \nu)$ decays too slowly to be integrable (that is no faster than $1/\nu$), its Fourier transform can no longer be understood as an integral and the identity~\eqref{eq:convTheorem} cannot be used to improve performance, since the FFT will fail to correctly determine either $\hat{f}(t)$ or $\hat{g}(t)$. Despite these complications, the convolution may still be defined as an ordinary integral (at least as long as the product $f( \nu)g( \nu)$ decays faster than $1/\nu$) and numerical evaluation is cumbersome but straightforward.\\
Here the LFT really excels. On the one hand, it can be evaluated much faster than any direct evaluation of the convolution (even if performed on an optimized grid) due to its superior scaling in the number of data points. On the other, it is able to cope with non-integrable functions as long as a suitable trade-off parameter exists. In other words, the LFT is much less plagued by convergence problems than the FFT. Slowly decaying functions, as we have already pointed out in the last section, can be numerically transformed with the LFT at the price of dynamical compression that inevitably reduces the signal-to-noise contrast at small arguments in the image space. However, for convolutions this problem is slightly less pronounced: If possible, setting the trade-off parameter of the final $\tau \to \omega$ back-transform in~\eqref{eq:convTheorem} to $k_\text{back}=1-k_1-k_2$,  where $k_1$ and $k_2$ are the optimized parameters for the transforms of $f(\eta e^{\omega})$ and $g(\eta e^{\omega})$, respectively, renders the result unaffected by any dynamical compression in $\tau$ because the noise level in the time domain stays constant at $\mathcal{O}{(\epsilon)}$, as the problematic exponential prefactors cancel. In the exotic case of strongly divergent integrals, this value of $k_\text{back}$ may not be useful if $e^{k_\text{back}\omega_\text{max}}$ becomes much larger than the expected result, in which case more data points and the less aggressive, symmetric choice $k_{1}=k_{2}=k_\text{back}$ are typically better suited (see last example in~\ref{sec:Ex}).

\section{Examples and optimizations}\label{sec:Ex}
\begin{figure}[tp]
\begin{center}
\includegraphics[width=0.9\columnwidth]{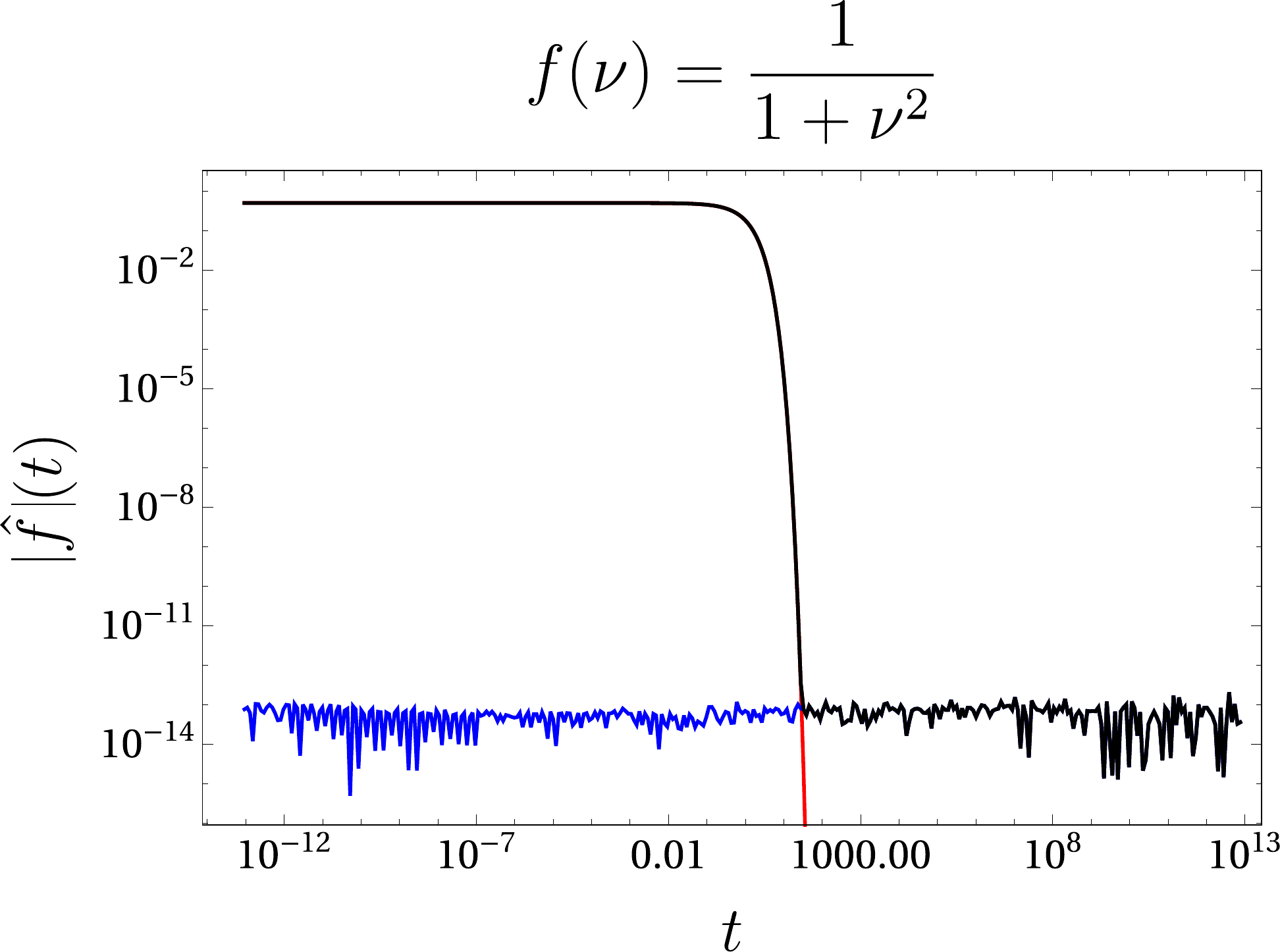}
\caption{(Color online) Fourier transform of $f(\nu)=\frac{1}{1+\nu^2}$ with $\Delta\omega=\Delta\tau=1/6$, $\Delta s=1/10$, $k=-1/100$ and symmetric intervals $\omega_s=s_s=\tau_s=-N/2$ on $N=360$ data points. The red line represents the analytical result $e^{-|t|}/(2\pi)$, while the numerical data is shown in black and the difference between the two in blue.}
\label{fig:1}
\end{center}
\end{figure}

\begin{figure}[tp]
\begin{center}
\includegraphics[width=0.9\columnwidth]{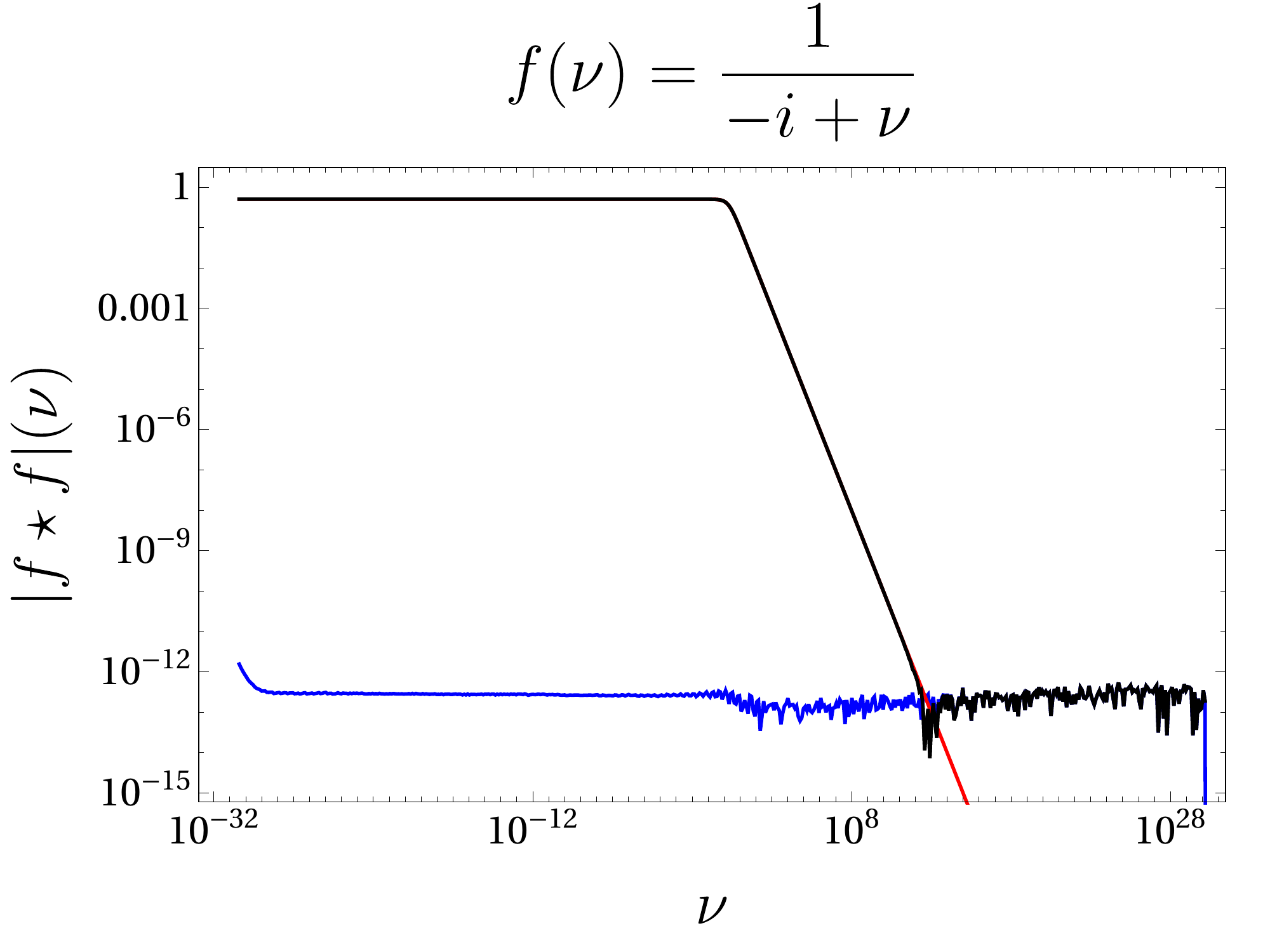}
\caption{(Color online) Convolution of $f(\nu)=\frac{1}{\nu-i}$ with itself, here $N=560$, $\Delta\omega=1/4$, $\Delta s=5/76$, $\Delta \tau=1/8$, $k_a=k_b=0.51$, $k_\text{back}=-0.02$, $\tau_s=-440$ and $s_s=\omega_s=-N/2$ were used. Color coding is the same as in Fig.~\ref{fig:1}.}
\label{fig:conv1}
\end{center}
\end{figure}

To benchmark the LFT and to illustrate the role of the transformation parameters, in particular of $k$,  we compute the Fourier transforms for several examples and compare the results to the exact solutions.
%Without loss of generality, we will assume that any shift in $\nu$ has already been removed, such that the functions we focus on are all centered around the origin and that $\nu$ furthermore has been rescaled via $\nu_0$ in such a way that the characteristic scale is unity.%
Without loss of generality, we focus on functions that are centered around the origin and that vary on a characteristic scale of unity. Deviations from that behavior can be remedied by preprocessing the function with a variable transformation which combines a shift of the original argument followed by rescaling it with a proper $\bar\nu$.\\
Let us begin with a benign example, a Lorentzian curve
\begin{align}\label{eq:ex1}
f(\nu)=\frac{1}{1+\nu^2}\,,
\end{align}
which could also be transformed with an ordinary fast Fourier transform. However, the slow convergence of the integral implies that reaching a global precision of 
$\epsilon = 10^{-12}$ with the FFT requires roughly $10^{13}$ data points, which exceeds numerical feasibility by several orders of magnitude. In contrast, to achieve the same accuracy with the LFT only a little more than 300 points suffice, as can be deduced from 
Eq.~\eqref{eq:minN}  with $R^{(1)}=\pi/2$. Indeed, Fig.~\ref{fig:1} has been obtained with 360 points and the (quasi-)optimal parameter $k=-1/100$, since $k_{\text{opt}}=0$ is prohibited by the $\Gamma$ function. Setting up the LFT in this way, the precision is no longer limited by the finite resolution, but by double-precision floating point arithmetic and error-propagation therein. In addition, the interval in $t$ can be chosen arbitrarily by adjusting $\Delta\tau$ and $\tau_s$ with no influence on the error level. As discussed in Sec.~\ref{sec:idealTOP}, due to the proximity to the pole of the $\Gamma$ function the last point $\hat{f}(t_\text{max})$ has to be subtracted which corresponds to the constant that arises from the leading $m=0$ contribution to $\exp(-k \tau)E_{2 \Gamma}^\sigma(\tau)$, see Eq.~\eqref{eq:E2Gamma}. The same procedure, which amounts to nothing else than a trivial subtraction of a one-parameter fit has been used for all other plots (except Fig.~\ref{fig:conv1}) as well.\\

\begin{figure}[tp]
\begin{center}
\includegraphics[width=0.9\columnwidth]{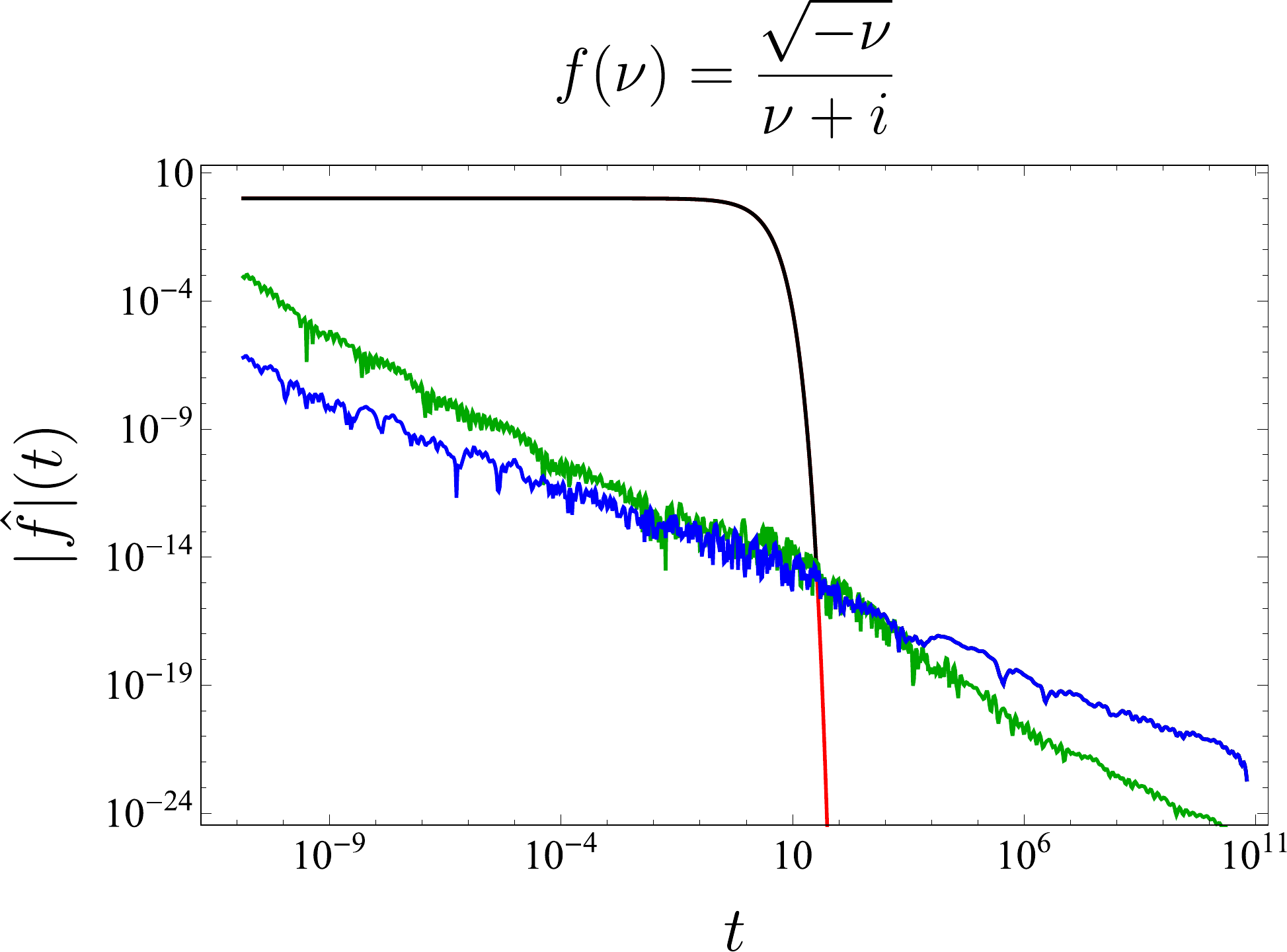}
\caption{(Color online) Numerical Fourier transform of $f(\nu)=\frac{\sqrt{-\nu}}{\nu+i}$ (black line) with $N=1000$ points, $\Delta\omega=1/5$, $\Delta s=2/45$, $\Delta\tau=1/20$, $s_s=\tau_s=-N/2$ and two different trade-off parameters. In green the difference between the exact analytical function $\hat{f}(t)$ (red line) and the numerical result with $k=1.01$ and $\omega_s=-N/2$ is shown, while the blue line depicts the same error but for the less aggressive $k=0.71$ and $\omega_s=-200$, which leads to smaller errors in the limit $t \to 0$ but to enhanced noise for $t \to \infty$, as discussed in the main text.} 
\label{fig:2}
\end{center}
\end{figure}

\begin{figure}[tp]
\begin{center}
\includegraphics[width=0.9\columnwidth]{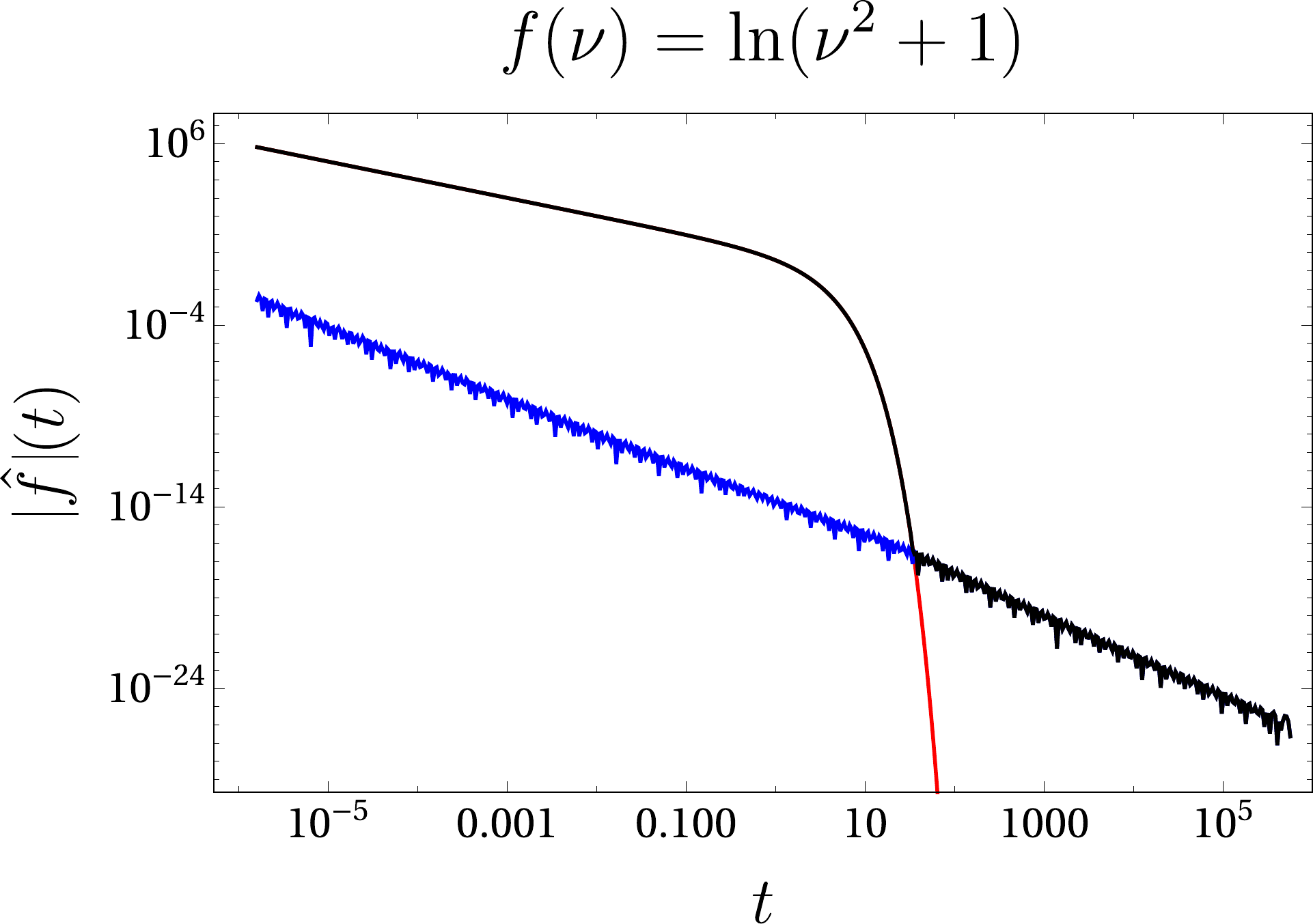}
\caption{(Color online) Fourier transform of $f(\nu)=\ln(\nu^2+1)$ with $N=560$, $\Delta\omega=1/7$, $\Delta s=1/14$, $\Delta \tau=1/21$, $k=2.05$ and $\omega_s=s_s=\tau_s=-N/2$. Color coding is the same as in Fig.~\ref{fig:1}.}
\label{fig:3}
\end{center}
\end{figure}

A significantly more demanding example (on the branch where $\sqrt{-1}=i$) is given by the function
\begin{align}
f(\nu)=\frac{\sqrt{-\nu}}{\nu+i}\,,
\end{align}
whose Fourier transform has to be understood in the sense of tempered distributions. For positive arguments $t>0$ it reads
\begin{align}
\hat{f}(t)=\frac{(1-i)}{\sqrt{2}}e^{-t}\,.
\end{align}
According to section~\ref{sec:idealTOP}, the optimized trade-off parameter is close to $k=1$, which removes the divergent behavior from the numerical integrals at the price of reduced precision at very small values of $t$. To demonstrate how the ideal choice of $k$ might depend on the data range of interest in the image space, Fig.~\ref{fig:2} depicts the result for  the ideal trade-off parameter $k=1.01$ and the suboptimal $k=0.71$. In order to achieve errors of $10^{-12}$ at $\tau=0$ in both cases, which requires roughly $N=600$ in the optimized setting, the grid size has been increased to $N=1000$. In agreement with the general discussion, values of $k$ larger than $1$ reduce errors at large arguments, while those smaller than unity increase precision close to $t=0$.\\
The next example shows
\begin{align}
f(\nu)=\ln(\nu^2+1)\,,
\end{align}
which transforms into
\begin{align}
\hat{f}(t)=\frac{e^{-|t|}}{|t|}
\end{align}
and is again correctly described by the LFT on only 560 points (see Fig.~\ref{fig:3}). However, the divergence of $f(\nu)$ as $\nu \to \infty$ results in an even stronger dynamical suppression than before. 

\begin{figure}[tp]
\begin{center}
\includegraphics[width=0.9\columnwidth]{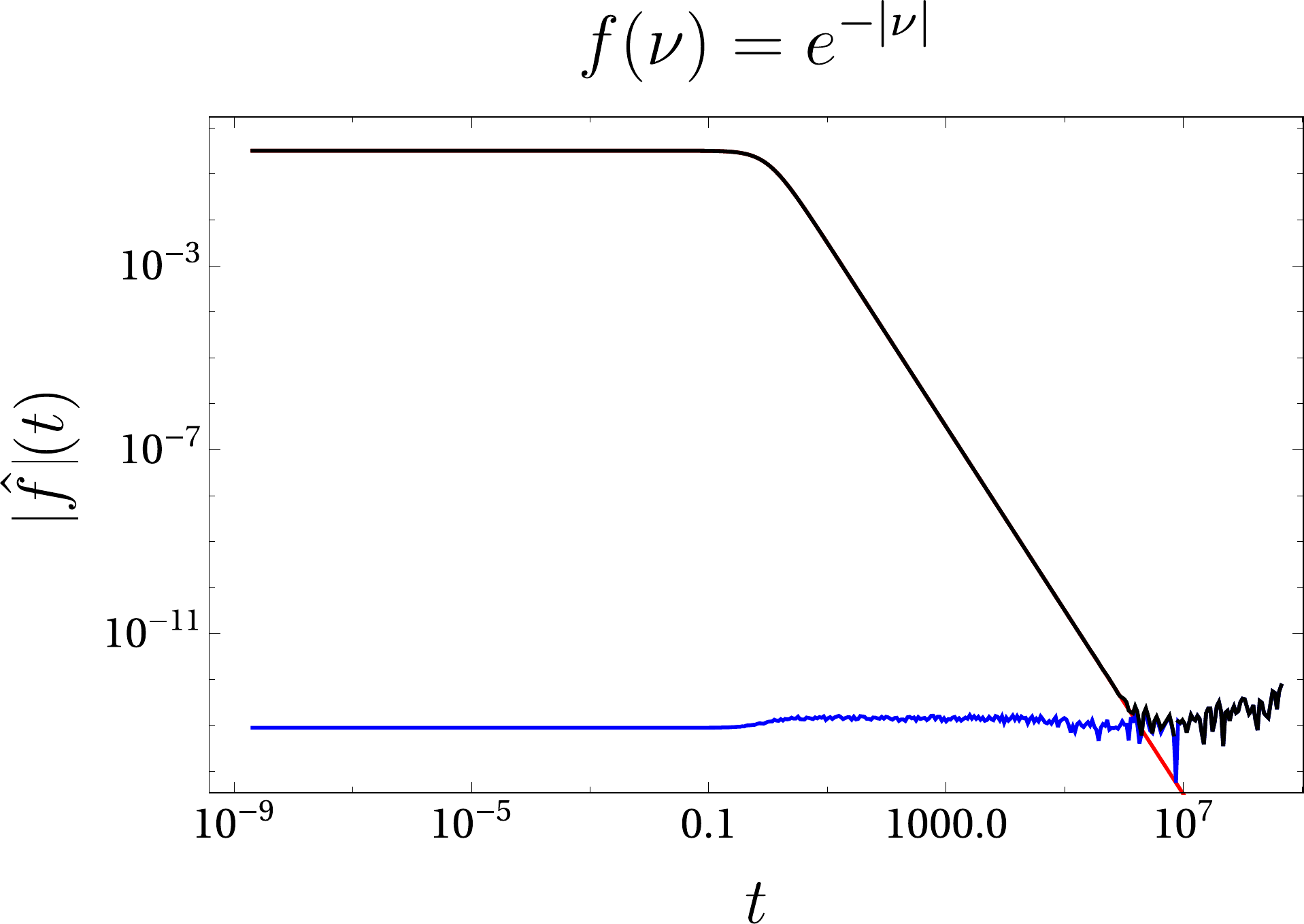}
\caption{(Color online) Fourier transform of $f(\nu)=e^{-|\nu|}$ with $N=480$, $\Delta\omega=1/15$, $\Delta s=2/21$, $\Delta \tau=1/12$, $k=-3/10$, $\omega_s=-420$ and $s_s=\tau_s=-N/2$. Color coding is the same as in Fig.~\ref{fig:1}.}
\label{fig:4}
\end{center}
\end{figure}

The well-behaved function $f(\nu)=e^{-|\nu|}$, which in fact corresponds to the inverse transformation of the very first example~\eqref{eq:ex1}, could also be treated by an FFT. However, the same accuracy as demonstrated in Fig.~\ref{fig:4} with $N=480$ on the logarithmic grid would require more than a million data points on a linear grid, illustrating that even for the most benign functions the LFT can outperform the direct application of an FFT.

%Finally, in Fig.~\ref{fig:5} we demonstrate, that also more structured functions like the Fourier transform of
%\begin{align}
%f(\nu)=\frac{1}{1+\nu^4}\, ,
%\end{align}
%which reads
%\begin{align}
%\hat{f}(t)=\frac{1}{4\sqrt{2}}e^{-\sqrt{i}|t|}\left((1+i)+(1-i)e^{i\sqrt{2}|t|}\right)\, ,
%\end{align}
%can be treated without loss of precision.\\
%\begin{figure}[htp]
%\begin{center}
%\includegraphics[width=0.9\columnwidth]{../plots/fig5.pdf}
%\caption{(Color online) Fourier transform of $f(\nu)=\frac{1}{1+\nu^4}$ with $N=600$, $\Delta\omega=1/15$, $\Delta s=2/15$, $\Delta \tau=1/12$, $k=-1/10$, $\omega_s=-435$, $\tau_s=-348$ and $s_s=-N/2$. Color coding is the same as in Fig.~\ref{fig:1}.}
%\label{fig:5}
%\end{center}
%\end{figure}
The convolution of 
\begin{align}
f(\nu)=\frac{1}{-i+\nu}
\end{align} 
with itself, that appears frequently in the evaluation of Feynman diagrams with non-relativistic propagators~\cite{abri75,fett71}, cannot be treated by FFTs, as the integral over $f(\nu)$ does not exist. Using the LFT remedies this issue by the help of the trade-off parameter. For instance, using the optimized value of $k$ given in section~\ref{sec:idealTOP}, a constant error of roughly $10^{-12}$, which is limited only by the internal floating point precision, is obtained with only $N=560$ points. In Fig.~\ref{fig:conv1} the two leading contributions to $E^\sigma_{2 \Gamma}$ from Eq.~\eqref{eq:E2Gamma} for the transformation from $t$ to $\nu$ with $m=0,1$ were subtracted by fitting the two corresponding parameters $I_1^\sigma(-i(k+m))$ to the high-frequency range near $\nu=10^{28}$.

\begin{figure}[tp]
\begin{center}
\includegraphics[width=0.9\columnwidth]{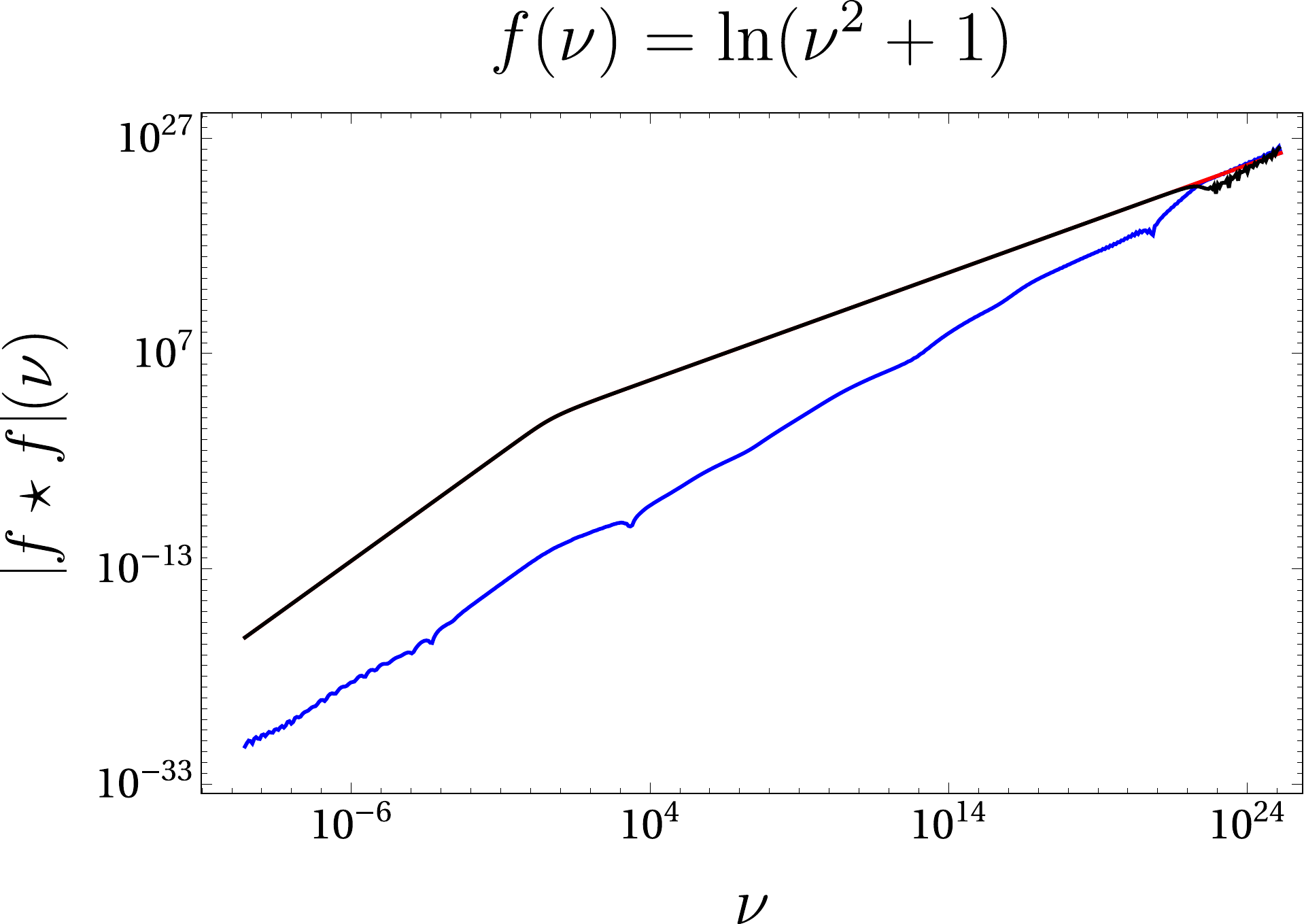}
\caption{(Color online) Convolution of $f(\nu)=\ln(\nu^2+1)$ with itself, here $N=560$, $\Delta\omega=1/7$, $\Delta s=1/14$, $\Delta \tau=1/14$, $k_a=k_b=k_\text{back}=8/5$, $\omega_s=-154$, $\tau_s=-495.6$ and $s_s=-N/2$ were used. Color coding is the same as in Fig.~\ref{fig:1}.}
\label{fig:conv2}
\end{center}
\end{figure}

Encouraged by these results, one can try to convolve some more exotic functions, for example $f(\nu)=\ln(\nu^2+1)$ with itself. In this case neither the convolution as an integral, nor the product of the distributions in the time domain is in general well-defined~\cite{GelfandBook}. Nevertheless, employing a cutoff $e^{-\delta |\nu|}$, with $\delta >0$, in the logarithmic frequency space and sending $\delta$ to zero at the end one finds analytically
\begin{align}
(f\star f)(\nu)=2\ln\left(1+\frac{\nu^2}{4}\right)-2\nu \arctan{\frac{\nu}{2}}\, .
\end{align}
As Fig.~\ref{fig:conv2} demonstrates, this result is again very accurately recovered by means of an LFT with $N=560$. Here, due to the large frequency interval used, round-off errors multiplied by $e^{-k_\text{back}\omega}$ are the limiting factor at large $\nu$. This affects the choice of trade-off parameters for the forward and backward LFTs, which is expected on general grounds, as remarked at the end of section~\ref{sec:Convo} on convolutions.
Furthermore, the first sub-leading divergence $\propto e^{2\omega}$ due to the next pole $(m=2)$ of the $\Gamma$ function beyond the constant has been fitted with a single parameter against the raw result at $\omega_N$ and subtracted.

\section{Application to physical examples}\label{sec:physEx}
Following the purely mathematical discussion, we now provide physical examples to highlight the real-world advantages of the LFT. These demonstrations are deliberately chosen to be simple, yet of relevance to current research and with apparent generalizations to more challenging problems.
\subsection{Polarization function}
The polarization function of the one-dimensional Bose gas is given by \cite{Mahan_book,Kamenev_book}
\begin{align}\label{eq:Lindhard}
\Pi(\omega,q)=-\int\frac{dk}{2\pi}\frac{n_B(\xi_k)-n_B(\xi_{k-q})}{\omega-\xi_k+\xi_{k-q}+i0^{+}}
\end{align}
with $n_B(k)=1/(e^{\beta \xi_k}-1)$ the Bose-Einstein distribution for the inverse temperature $\beta=1/(k_B T)=1$ and the free dispersion $\xi_k=k^2/(2m)-\mu$ with momentum $k$, mass $m=1/2$ and chemical potential $\mu$. $\Pi(\omega,q)$ describes density-density correlations and in general has no known closed expression. One therefore has to rely on a numerical evaluation of the integral. In case of the fugacity $z=\exp{(\beta\mu)}$ approaching unity from below, the density fluctuations in the regime of long wave lengths proliferate, which is reflected in the singular behavior $n_B(0)\sim-1/(\beta\mu)$. As a consequence, the direct evaluation of the polarization function becomes numerically expensive. However, precisely these low temperature correlation functions are a common ingredient in quantum many-body theories, in particular: quantum critical transport~\cite{sach11}, response near phase transitions in ultracold atoms~\cite{endr12} and Bose gases in optical cavities~\cite{rits13}.

\begin{figure}[tp]
\begin{center}
\includegraphics[width=0.9\columnwidth]{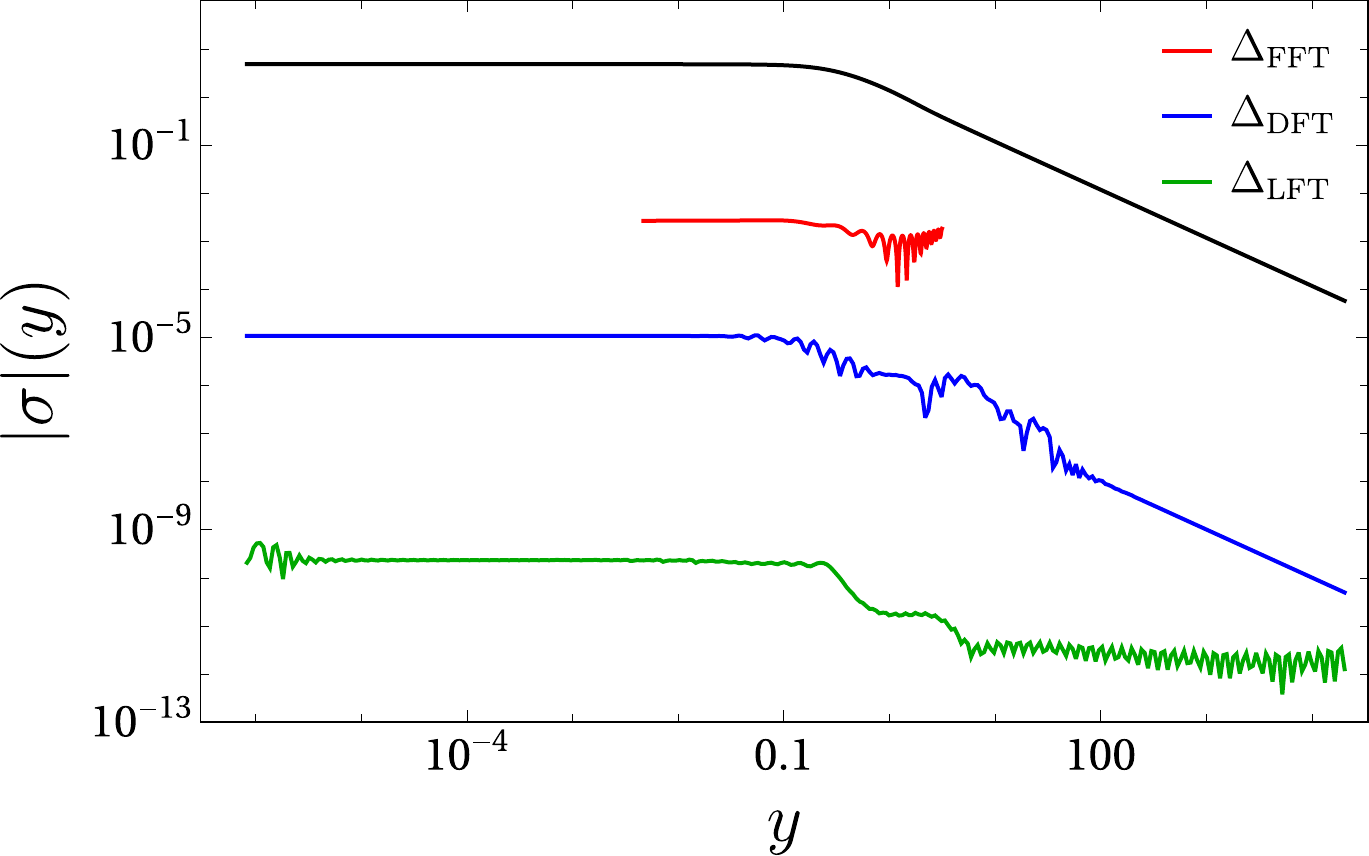}
\caption{(Color online) Comparison of the absolute error in the evaluation of $\sigma(y)$ with $\mu=-1/10$ using an FFT (red) a DFT following a cubic spline interpolation (blue) and the LFT (green). For comparison we also show $\sigma(y)$ in black. The parameters of the LFT are $N=336$, $\Delta\omega=1/14$, $\Delta s=1/4$, $\Delta \tau=1/11$, $k=-9/10$, $\omega_s=-196$, $\tau_s=-234$ and $s_s=-N/2$ for the first and the same values, except $\Delta s=1/6$ and $k=1/20$, for the second transform. The same grid is then also used for the DFT. The FFT is run on a grid with the same number of points and a lattice spacing $\delta x\approx 0.07$ and $\delta y\approx 0.01$ for $x$ and $y$ in Eq.~\eqref{eq:Lindhard_Fourier}.}
\label{fig:lindhard}
\end{center}
\end{figure}

In the following, we will show that an efficient evaluation of $\Pi(\omega,q)$ with unrivaled precision is possible by means of the LFT. The integral in Eq.~\eqref{eq:Lindhard} can be rewritten as a convolution, which allows for an efficient treatment that requires only two one-dimensional (half-sided) Fourier transforms:
\begin{subequations}
\begin{align}
\!\Pi(\omega,q)&=\sum_{\eta=\pm}\frac{i\eta}{2q}\sigma\!\left(\frac{\omega+\eta q^2}{2q}\right)\\
\sigma(y)&=\mathcal{F}^{-1}_{x\to y}\left[\theta(-x)\mathcal{F}_{k\to x}\left(n_B(k)\right)(x)\right](y)\,.\label{eq:Lindhard_Fourier}
\end{align}
\end{subequations}
The computation of the polarization function on a two-dimensional $(\omega,q)$-grid of size $N\times N$ therefore requires only $\mathcal{O}(N^2)$ operations for the evaluation of $\sigma(y(\omega,q))$ in contrast to $\mathcal{O}(N^3)$ for the direct approach. In fact, the actual fast Fourier transform in Eq.~\eqref{eq:Lindhard_Fourier} results only in subleading corrections to the overall complexity. The main advantage of the LFT over other Fourier transforms, however, lies in its accuracy. We highlight this in Fig.~\ref{fig:lindhard} by comparing the absolute error obtained for an ordinary FFT, the LFT and a discrete Fourier transform (DFT) in combination with a cubic spline interpolation. All algorithms use the same number of data points as well as optimized transformation parameters, with the DFT and LFT operating on a common grid. For all relevant values the LFT outperforms the other methods by several orders of magnitude, furthermore, in contrast to the FFT a much larger interval can be sampled.

\subsection{Glass transition}
The glassy, mechanically rigid state of amorphous materials and its realization by supercooling liquids has been studied for a very long time~\cite{debe01}. Nevertheless, a theoretical description of this state is difficult since one has to deal with density fluctuations on various length scales and very slow relaxation processes, as observed in experiments~\cite{stef94}. One theoretical approach to this problem is mode coupling theory~\cite{reic05,goet08} which is based on an effective equation of motion for the dynamical structure factor $S(\mathbf k ,t)$. In the following, we illustrate how the LFT, which by construction is capable of dealing with multi-scale problems, can be applied to these kinds of models. Here we focus on a simplified variant of mode-coupling theory due to Leutheusser~\cite{leut84} and Bengtzelius et al.~\cite{beng84} and remark on the advantageous properties of the LFT for more generic problems of this kind.

To study the relaxation of density distortions one introduces an effective temporal order-parameter $\Phi(t)$, whose phenomenological time evolution is given by~\cite{leut84}
%\begin{align}
%\label{eq:EOM}
%\ddot{\Phi}(t) + \gamma \dot{\Phi}(t) + \Omega_0^2 \Phi(t) + 4 \Omega_0^2 \lambda \int_0^{t} d\tau \; %\Phi^2(\tau) \dot{\Phi}(\tau- t) = 0 \, .
%\end{align} 
\begin{align}
\label{eq:EOM}
\ddot{\Phi}(t) + \gamma \dot{\Phi}(t) + \Omega_0^2 \Phi(t) =- 4 \Omega_0^2 \lambda \int_0^{t} d\tau \; \Phi^2(\tau) \dot{\Phi}(\tau- t) \, .
\end{align}
The left side of the equation is a simple harmonic oscillator with damping rate $\gamma$ and frequency $\Omega_0$. The correlations responsible for the glass transition are incorporated in the memory integral on the right-hand side and weighted by the dimensionless, positive coupling constant $\lambda$. The initial conditions $\Phi(0)=1$ and $\dot{\Phi}(0)=0$ model the original deviation from the equilibrium state $\Phi \equiv 0$. Despite its simplicity, the above equation takes the basic properties of the glass transition into account which we briefly review before presenting the solution based on the LFT.

The physical order parameter distinguishes between two phases via the long-time limit $\Phi(t \to \infty) =c$: The ergodic phase is characterized by perfect relaxation corresponding to $c=0$, in contrast to the glass phase where the initial distortion never disappears completely and thus $c>0$ . To gain further insight into the phase diagram we apply the half-sided Fourier transformation 
\begin{align}\label{eq:halfFT}
\hat\Phi(\nu) = \int_0^\infty dt\,\Phi(t) e^{i  \nu t} \, . 
\end{align}
Furthermore, we reparametrize the function $\Phi(t) = \delta \Phi(t) +c$ to obtain the asymptotics
\begin{align}\label{eq:ansatz}
\delta \Phi(t) \to 
\begin{cases}
1-c & t\to 0^+ \\
0 & t \to \infty
\end{cases}\, ,
\end{align}
irrespective of the phase.
The constant value $c$ gives rise to a term proportional to $\delta(\nu)$ in the Fourier transform of  Eq.~\eqref{eq:EOM} that has to be canceled to satisfy $\delta \Phi(t \to \infty)=0$. This is the case if
\begin{align}
c(\lambda) = 
\begin{cases}
0  & \lambda < 1 \\
 \cfrac{1+\sqrt{1-1/\lambda}}{2} & \lambda \geq 1
 \end{cases}\,,
\end{align}
which not only determines the asymptotic value of the order parameter but also identifies the critical coupling for the glass transition $\lambda_c=1$. As has been shown in Refs.~\cite{leut84, beng84} by analytic means, the approach to the phase boundary is characterized by a divergent low-frequency limit of $\delta \hat\Phi(\nu)$ that follows the power law 
\begin{align}\label{eq:critExp}
\delta \hat\Phi(\nu = 0) \sim 
\begin{cases}
(1-\lambda)^\mu\!\! & \text{for } \lambda \to 1^- \,\text{with } \mu = 1.76498...\\
(\lambda-1)^{\mu^\prime}\!\!\! & \text{for } \lambda \to 1^+\, \text{with } \mu^\prime = 0.76498...
\end{cases}
\end{align}
In the ergodic phase the half-sided Fourier transform of Eq.~\eqref{eq:EOM} yields the self-consistent relation
\begin{align} \label{eq:EOMnormal}
\delta\hat\Phi(\nu) = - \frac{1}{i \nu + \displaystyle \frac{\Omega_0^2}{i \nu - \gamma - 4 \lambda \Omega_0^2 \mathcal{F}_{t\to\nu}[\delta \Phi^2(t)](\nu)}} \, .
\end{align}
Similarly, in the glass phase one obtains a quadratic equation with the solution
\begin{align}\label{eq:EOMglass}
\begin{split}
\delta \hat\Phi(\nu)& = \frac{- B \pm\sqrt{B^2 -4 A C}}{2 A}\\
A & = 8 c \lambda \Omega^2_0 i \nu \\
B & = \nu^2 +i \gamma \nu +\Omega_0^2\left[1-4 \lambda c^2 +  4 \lambda i \nu \mathcal{F}_{t\to\nu}[\delta \Phi^2(t)](\nu)\right]  \\
C & = (1-c) \left[-i \nu + \gamma +4 \lambda \Omega_0^2 \mathcal{F}_{t\to\nu}[\delta \Phi^2(t)](\nu)\right]  \, ,
\end{split}
\end{align} 
where in the first line one has to choose the branch that yields a positive $\Re \hat\Phi(\nu)$, since this function represents a retarded, bosonic correlation function~\cite{fett71}. The equations~\eqref{eq:EOMnormal} or~\eqref{eq:EOMglass} can be solved in an iterative manner, which requires repeated Fourier transformations between the time and frequency spaces. To reliably compute the critical behavior in the vicinity of the glass transition, however, one has to include very long times, especially to capture the extremely slow relaxation when approaching the instability from the ergodic phase. Yet, this is exactly the scenario the LFT has been devised for.

\begin{figure}[tp]
\begin{center}
\includegraphics[width=0.9\columnwidth]{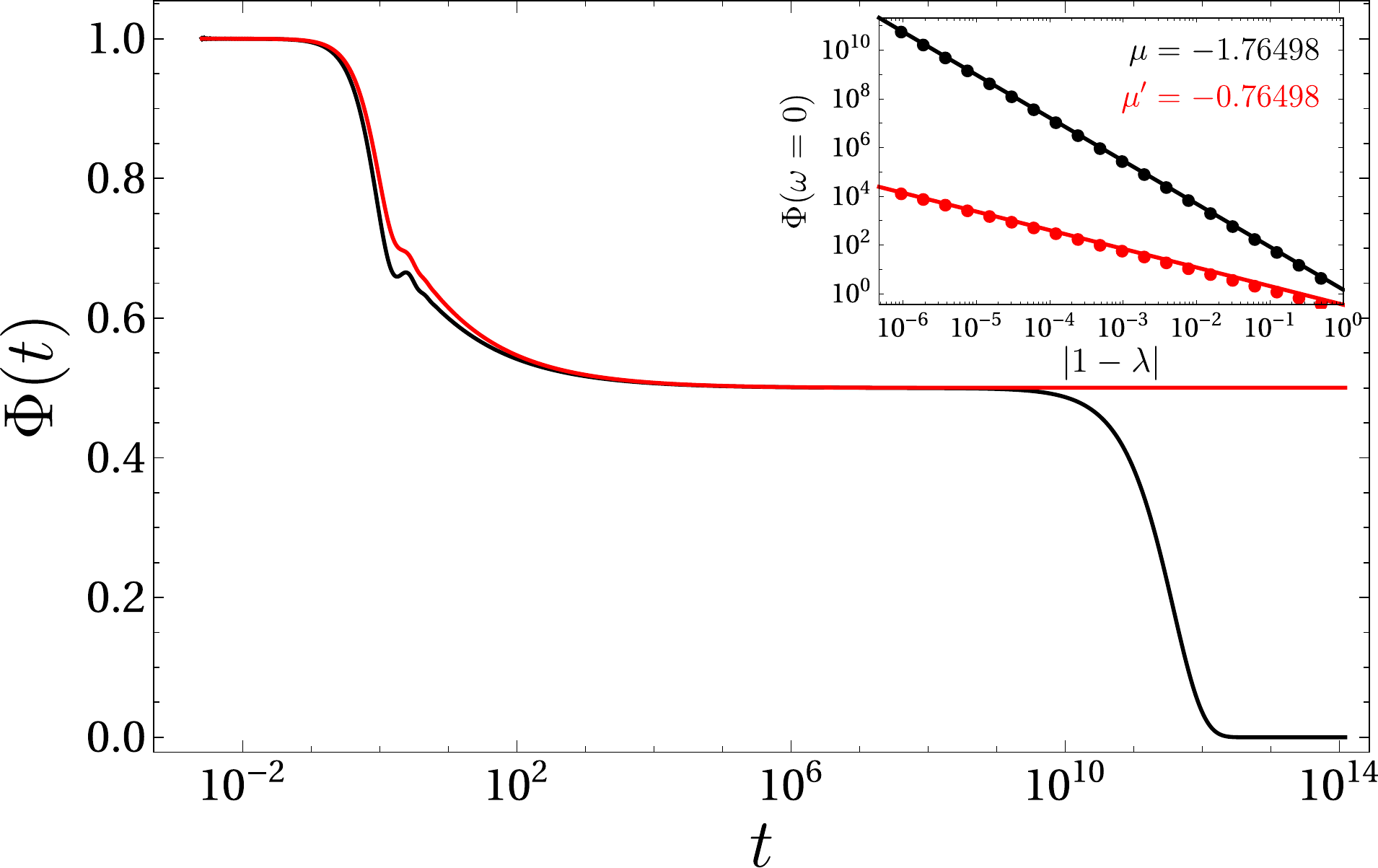}
\caption{(Color online) Time dependence of the order parameter for $\lambda = 1-2^{-21}$ (black) in the ergodic phase and for $\lambda=\lambda_c+2^{-21}$ (red) in the glass phase at $\gamma=1 = \Omega_0$. Inset: data points represent the numerical results for $\delta\hat\Phi(0)$. The straight lines are the analytical results for the critical exponent below (black) and above (red) the critical point. All results have been obtained with $N=2400$.}
\label{fig:glass}
\end{center}
\end{figure}

Note that the formulation of Eq.~\eqref{eq:EOM} in terms of $\delta \Phi$ is not only helpful for analyzing the problem in further detail, but in order to apply the LFT it is also mandatory because $\delta \Phi$ satisfies the condition~\eqref{eq:condition}. 
 Figure~\ref{fig:glass} shows $\Phi(t)$ at the values $\lambda =1 \pm 2^{-21}\approx 1\pm 4.7\cdot 10^{-7}$ in the immediate vicinity of the glass transition. In the ergodic phase the plateau, which characterizes the so-called regime of $\beta$-relaxation, reaches times of $\mathcal O(10^{12})$ before the final $\alpha$-relaxation to $\Phi=0$ sets in. In addition to the global features of the dynamics, the LFT reproduces the critical exponents from Eq.~\eqref{eq:critExp} with a numerical error on the order of $10^{-4}$ (see inset of Fig.~\ref{fig:glass}). 
To achieve such small errors one can profoundly benefit from the flexibility of the LFT: While 
the LFTs have to operate on identical $\omega$ and $\tau$ grids, irrespective of the direction of the transformation, the auxiliary $s$ space can be sampled for the two directions $\omega \to \tau$ and $\tau \to \omega$ independently. Moreover, one can introduce two LFTs for $\Re \delta \Phi(\nu)$ and $\Im \delta \Phi(\nu)$ separately without altering the overall computational cost. 
%\jls{We call the resulting shifts $s_{R,I,\tau}$ and the step widths $\Delta s_{R,I,\tau}$. In the same spirit, one independently defines the trade-off parameters $k_R, k_I$ and $k_\tau$  for the three LFTs, which is quite helpful to optimally capture the individual asymptotics of the three functions. Secondly, for each value of $\lambda$ one uses the shifts $\omega_s$ and $\tau_s$ to minimize the truncation errors. }
In total, the numerical effort to obtain $\Phi(t)$ at all times for a given $\lambda$ scales like $N_{\text{it}} \cdot N \log N$  where $N_{\text{it}}$ denotes the number of iterations needed to reach convergence.
Using optimized trade-off parameters and a grid of length $N=2400$ suffices to produce the results shown in Fig.~\ref{fig:glass}. The small number of data points used in the LFT allows to find the converged order parameter in a couple of seconds, even close to the phase transition, where $N_{\text{it}} \sim 10^4$.
A direct numerical solution of the integro-differential equation~\eqref{eq:EOM} on a discretized time axis requires step sizes $\Delta t \lesssim 1$ independent of the magnitude of $t$ in order to compute the cancellations between the various terms with sufficient precision. These are responsible for the slow evolution and the formation of the plateau over six orders of magnitude in time. Larger $\Delta t$ would lead to an  instability of the numerical solution that diverges away from the physical $\Phi(t)$. Due to the scaling $\mathcal O (N_{\text{dir}}^2)$ of a direct approach and limited step size it is numerically completely unfeasible to reach times of order $10^{14}$. This, however, is necessary to reliably determine $\Phi(\omega = 0)$ and consequently the critical exponents \footnote{It is worth to mention, that in combination with a partially analytical solution, step sizes $\Delta t \gg 1$ are possible. The implementation of such a method is, however, far more complicated and time-consuming than the simple, direct approach presented here~\cite{Goetze1988}.
%If I understand Haussmann correctly he solves the equation with a fixed integral and then repeatedly improves the estimate for the next step before moving on. He also uses an adaptive grid and repeated spline interpolations.
}.

Regarding more complicated versions of mode-coupling theory that resolve the dependence on the length scales, thereby considering the structure factor $S(\mathbf{k}, t)$ instead of $\Phi(t)$, the grid size of $N \sim 10^3$ used here is still small enough to incorporate a second argument without running into memory limitations. If the coupling between different wave vectors can be written in terms of convolutions, which usually is the case~\cite{reic05,goet08}, the LFT can also be applied to simplify the spatial dependence. Very similar problems appear in the context of approximate equations of motions of correlation functions in quantum field theory, which typically include algebraic decays in frequency and momentum space. As mentioned earlier, an example of the application of the LFT in the context of ultracold Fermi gases can be found in Ref.~\cite{Frank2018}.

\section{Conlcusion}\label{sec:con}
We have shown rigorously that the LFT can be used to numerically transform nonintegrable functions, as long as their asymptotics can be controlled by the trade-off parameter. Furthermore, we have proven that one can achieve exponential convergence in the number of data points if the function is analytic in a cone with finite opening angle around the real axis in the original argument. Finally, we have given several examples that benchmark the superior convergence of the LFT compared to the FFT including functions that have to be considered within the concept of generalized Fourier transformations. \\
\emph{Acknowledgments} The authors thank Wilhelm Zwerger for fruitful discussions and comments on the manuscript and Christian Johansen for pointing out erroneous factors of $2\pi$. This work has been supported by the Nanosystems Initiative Munich (NIM).
\bibliographystyle{unsrt}
\bibliography{biblio_v2}

%Certain Meijer G functions can be calculated very efficiently with the LFT to arbitrary precision: Using the convolution theorem for Meijer G functions (see Prudnikov 1990: Integrals and Series Volume 3 p. 350) complicated functions can be formulated as the Fourier transform of simple G functions. Repeating the procedure an entire stack of Meijer G functions can be calculated. Typically the direct evaluation of the integral form of these functions is very time-consuming, often rendering the Meijer G function a dead end for practical quasi-exact calculations. Due to the excellent convergence properties of the LFT rapidly evaluated representations for many G functions can be constructed. Increasing the internal floating point precision there is no limit to the accuracy of this approach and the fast evaluation of each transform allows to stack many of them to reach even the most unpleasant among these special functions. Actually during the numerical tests for this publication the built-in evaluation of \emph{Mathematica} for the Meijer G function was almost always slower and less precise than the corresponding LFT calculation.

\end{document}